\documentclass[10pt,english]{article}
\usepackage[T1]{fontenc}
\usepackage[latin1]{inputenc}
\usepackage{a4wide}
\usepackage{babel}
\usepackage{latexsym}
\usepackage{color}
\usepackage{hyperref}
\usepackage{amsthm}

\makeatletter

\newtheorem{prop}{Proposition}[section]
\newtheorem{lemma}[prop]{Lemma}
\newtheorem*{lem}{Lemma}
\newtheorem{theorem}[prop]{Theorem}
\newtheorem{defi}[prop]{Definition}
\newtheorem{rem}[prop]{\em Remark\/}
\newtheorem{cor}[prop]{Corollary}
\newtheorem{cond}[prop]{Conditions}
\newtheorem{ex}[prop]{\em Example\/}
\newtheorem{sch}[prop]{Scholium}

\newcommand{\bp}{\begin{prop}}
\newcommand{\bl}{\begin{lemma}}
\newcommand{\bt}{\begin{theorem}}
\newcommand{\bd}{\begin{defi}\rm}
\newcommand{\br}{\begin{rem}\rm}
\newcommand{\be}{\begin{equation}}
\newcommand{\bea}{\begin{eqnarray}}
\newcommand{\bpr}{\begin{proof}}
\newcommand{\bc}{\begin{cor}}
\newcommand{\bco}{\begin{cond}}
\newcommand{\bex}{\begin{ex}\rm}
\newcommand{\bsch}{\begin{sch}}

\newcommand{\ep}{\end{prop}}
\newcommand{\el}{\end{lemma}}
\newcommand{\et}{\end{theorem}}
\newcommand{\ed}{\end{defi}}
\newcommand{\er}{\end{rem}}
\newcommand{\ee}{\end{equation}}
\newcommand{\eea}{\end{eqnarray}}
\newcommand{\epr}{\end{proof}}
\newcommand{\ec}{\end{cor}}
\newcommand{\eco}{\end{cond}}
\newcommand{\eex}{\end{ex}}
\newcommand{\esch}{\end{sch}}

\newcommand{\nn}{ \nonumber \\ }
\newcommand{\pr}{{\em Proof.\ }}
\renewcommand{\qed}{\vrule height 5pt width 5pt depth 0pt}
\newcommand{\ui}{^{(1)}}
\newcommand{\uii}{^{(2)}}
\newcommand{\di}{_{(1)}}
\newcommand{\dii}{{}_{(2)}}
\newcommand{\ot}{\otimes}

\newcommand{\asrd}{{{\cal A}^*}}
\newcommand{\atrd}{{^*\!\! {\cal A}}}
\newcommand{\atld}{{{\cal A}_*}}
\newcommand{\asld}{{_*{\cal A}}}
\newcommand{\ci}{\circ}
\newcommand{\lu}{\leftharpoonup}
\newcommand{\ld}{\leftharpoondown}
\newcommand{\ru}{\rightharpoonup}
\newcommand{\rd}{\rightharpoondown}

\newcommand{\End}{{\rm End}} 
\newcommand{\Hom}{{\rm Hom}} 
\newcommand{\M}{\mathcal{M}}
\newcommand{\A}{\mathcal{A}}
\newcommand{\stac}[1]{\stackrel{\ot}{_{_{#1}}}}
\newcommand{\op}{^{op}}
\newcommand{\pri}{^{\prime}}
\newcommand{\inv}{^{-1}}
\newcommand{\h}{\ {\bf\underline{\ }}\ }

\newcommand{\il}{{\cal L}}
\newcommand{\ir}{{\cal R}}
\newcommand{\err}{\wp}

\newcommand{\setc}[1]{\setcounter{equation}{#1}}

\newcommand{\lb}{\label}

\begin{document}

\large
\title{\bf Integral theory for Hopf Algebroids}
 
\author{\sc Gabriella B\"ohm \\
Research Institute for Particle and Nuclear Physics, Budapest,\\
H-1525 Budapest 114, P.O.B. 49, Hungary\\
E-mail: G.Bohm@rmki.kfki.hu}

\date{}
 
\maketitle
\normalsize 
\begin{abstract}  
\noindent
The theory of integrals is used to analyse the structure of Hopf algebroids
\cite{B,HGD}. We prove that the total algebra of a Hopf algebroid is a
separable extension of the base algebra if and only if it is a
semi-simple extension and  if and only if the Hopf algebroid
possesses a normalized integral. 
{
The total algebra of a finitely generated and projective Hopf algebroid 
is a Frobenius extension of the base algebra} 
if and only if
the Hopf algebroid possesses a non-degenerate integral. We give also a
sufficient and necessary condition in terms of integrals, under which it is a
quasi-Frobenius extension, and illustrate by an example that this
condition does not hold true in general. Our results are
generalizations of classical results on Hopf algebras \cite{LaSwe,Par}. 
\end{abstract}
\bigskip
\bigskip
\bigskip
\bigskip

\section{Introduction}

The notion of {\em integrals} in Hopf algebras
has been introduced by Sweedler \cite{Swe69}. The integrals in Hopf
algebras over principal ideal domains were analysed in
\cite{LaSwe,Swe}  where the 
following -- by now classical -- results have been proven: 

-- A free, finite dimensional bialgebra over a principal ideal domain
   is a Hopf algebra if and only if it possesses a non-degenerate left
   integral.
   (Larson-Sweedler Theorem.)

-- The antipode of a free, finite dimensional Hopf algebra over a principal
   ideal domain is bijective.

-- A Hopf algebra over  a field is finite dimensional if and only if
  it possesses a non-zero left integral.

-- The left integrals in a finite dimensional Hopf algebra over a field
  form a one dimensional subspace.

-- A Hopf algebra over a field 
  is semi-simple if and only if it possesses a normalized left integral.
  (Maschke Theorem.)

There are numerous generalizations of these results in the
literature. Historically the first is due to Pareigis \cite{Par} who proved
the following statements on a finitely generated and projective Hopf algebra
$(H,\Delta,\epsilon,S)$ over a commutative ring $k$: 

-- $H$ is a Frobenius extension of $k$ if and only if 
   there exists a Frobenius functional $\psi:H\to k$ satisfying
   $(H\ot \psi)\ci \Delta=1_H \psi(\h)$.

-- The antipode, $S$, is bijective. 

-- The left integrals form a projective rank 1 direct summand of the
   $k$-module $H$. 

-- $H$ is a quasi-Frobenius extension of $k$. 

-- A  finitely generated and projective bialgebra  over a commutative ring $k$,
   such that 
   $\textrm{pic}(k)=0$, is a Hopf algebra if and only if it possesses a
   non-degenerate left integral.

The generalization of the Maschke theorem to Hopf algebras $H$ over
commutative rings $k$  states that the existence of a normalized left
integral in $H$ is equivalent to the separability of $H$ over $k$,
what is further equivalent to its relative semi-simplicity in the
sense \cite{Hattori,HirSug} that any $H$-module is $(H,k)$-projective
\cite{CaenMi,Lomp}. This is equivalent to the true 
semi-simplicity of $H$ (i.e. the true projectivity of
any $H$-module \cite{Pierce}) if and only if $k$ is a semi-simple ring
\cite{Lomp}. 

As a nice review on these results we recommend  Section 3.2 in
\cite{CaMiZhu}. 

Similar results are known also for generalizations of Hopf
algebras. Integrals for finite dimensional quasi-Hopf algebras \cite{Dri} over
fields were studied in \cite{HN,Pan,PanVO,BuCae} and for finite dimensional
weak Hopf algebras \cite{BSz,BNSz} over fields in \cite{BNSz,Vecsi}.

The purpose of the present paper is to investigate which of the above results
generalizes to {\em Hopf algebroids}.

Hopf algebroids with bijective antipode have been introduced in
\cite{HGD,B}. It is important to emphasize that this notion of a Hopf algebroid
is not equivalent to that introduced under the same name by Lu in
\cite{Lu}. Here we generalize the definition of \cite{HGD,B} by relaxing the
requirement of the bijectivity of the antipode. A Hopf algebroid consists of
a compatible pair of a left and a right bialgebroid structure
\cite{T,Lu,Sz,Sz2} on the 
common total algebra $A$. The antipode relates these two left- and righthanded
structures. Left/right integrals {\em in} a Hopf algebroid are defined
as the invariants of the left/right regular $A$-module in terms of the
counit of the left/right bialgebroid. Integrals {\em on} a Hopf
algebroid are the comodule maps from the total algebra to the base algebra
(reproducing the integrals {\em in} the dual bialgebroids, provided
the duals possess bialgebroid structures).

The total algebra of a bialgebroid can be looked at as an extension of
the base 
algebra or its opposite via the source and target maps,
respectively. In this way there are four algebra extensions associated to
a Hopf algebroid. The main results of the paper relate properties of these
extensions to the existence of integrals with special properties: 
\begin{itemize}
\item A Maschke type theorem, proving that the separability, and also the
(in two cases left in two cases right)  
semi-simplicity of any of the four extensions
is equivalent to the existence of a normalized integral {\em in} the Hopf
algebroid (Theorem \ref{maschkethm}).
\item 
{ 
The total algebra is a Frobenius extension both of the base algebra 
and of its opposite algebra if and only if there exists a non-degenerate left
(equivalently, right) integral {\em in} the Hopf algebroid (Corollary
\ref{cor:Frob}). In particular, if the total algebra is a finitely generated
  and projective module of the base algebra (in several appropriate senses),
  then}
any of the four extensions is a Frobenius extension if and only
if there exists a non-degenerate 
{ (left or right)}
integral {\em in} the Hopf algebroid (Theorem \ref{frobthm}).
\item 
{
Under the same finitely generated projectivity conditions in the previous item,}
any of the four extensions is (in two cases a left in two cases a
right) quasi-Frobenius extension 
  if and only if 
  the (left or right) integrals {\em on} the Hopf algebroid form
  a flat module over the  base algebra
 (Theorem \ref{rightqf}).
\end{itemize}

Our main tool in proving the latter two points is the Fundamental Theorem for
Hopf modules over Hopf algebroids (Theorem \ref{fundi}).
\footnote{
{
In contrast to the Fundamental Theorem for Hopf modules over Hopf algebras,
Theorem \ref{fundi} is {\em not} proven for {\em arbitrary} Hopf algebroids. A
weaker version of the theorem, relying on a more restrictive notion of a
comodule of a Hopf algebroid, is proven for an arbitrary Hopf algebroid in
Theorem 3.26 and Remark 3.27 of the arXiv version of \cite{Bohm:rev}.
}}

The paper is organized as follows: We start Section \ref{int} with
reviewing some results on bialgebroids from
\cite{T,Lu,Sz,Sz2,Sch,BrzeMi,KSz,Sch1,Sz3}, the knowledge of which is
needed for the understanding of the paper. Then we present the definition
of Hopf algebroids and discuss some of its immediate
consequences. Integrals both {\em in} and {\em on} Hopf algebroids are
introduced and some equivalent characterizations are given.

In Section \ref{maschke} we prove two Maschke type theorems. The first one
collects some equivalent properties (in particular separability)
of the inclusion of the base algebra in the total algebra of a Hopf
algebroid. These equivalent properties are related to the 
existence of a normalized integral {\em in} the Hopf algebroid. The
second theorem collects some equivalent properties (in particular
coseparability) 
of the coring underlying a Hopf algebroid. These equivalent
properties are shown to be equivalent to the 
existence of a normalized integral {\em on} the Hopf algebroid. 

In Section \ref{frob} we prove the Fundamental Theorem for Hopf
modules over a Hopf algebroid. This theorem is somewhat stronger than
the one that can be obtained by the application of (\cite{Brze},
Theorem 5.6) to the present situation. 
{
Still, it is not known to hold for an arbitrary Hopf algebroid.
}
The main result of the section
is Theorem \ref{frobthm}. In proving it we follow an analogous
line of reasoning as in \cite{LaSwe}. That is, assuming that
{ 
all the four}
module structures of the total algebra over the base algebra 
{
are}
finitely generated and projective, we apply the Fundamental Theorem to
the Hopf module, constructed on the dual of the Hopf algebroid
(w.r.t. the base algebra). Similarly to the case of Hopf algebras, our
result implies the existence of non-zero integrals {\em on} any
finitely generated projective Hopf algebroid. Since the dual of a
(finitely generated projective) Hopf algebroid is not known to be a
Hopf algebroid in general, we have no dual result, that is, we do
not know whether there exist non-zero integrals {\em in} any finitely
generated projective Hopf algebroid. As a byproduct, also a sufficient
and necessary condition on a finitely generated projective Hopf
algebroid is obtained, under which the antipode is bijective. We do not
know, however, whether this condition follows from the axioms.

In Section \ref{qf} we 
use the results of Section \ref{frob} to obtain conditions which are
equivalent to the (either left or right) quasi-Frobenius property of
any of the four extensions  behind a 
{
finitely generated and projective}
Hopf algebroid.
In order to show that these conditions do not hold true in general, we
construct a counterexample.

Throughout the paper we work over a commutative ring $k$. That is, the
total and base algebras of our Hopf algebroids are $k$-algebras. For an
(always associative and unital) $k$-algebra $A\equiv(A,m_A,1_A)$ we
denote by $_A\M$, 
$\M_A$ and ${_A\M_A}$ the categories of left, right, and bimodules over
$A$, respectively. For the $k$-module of morphisms in $_A\M$,
$\M_A$ and ${_A\M_A}$ we write $_A\Hom(\ ,\ )$, $\Hom_A(\ ,\ )$ and
${_A\Hom_A}(\ ,\ )$, respectively. 

{\bf Acknowledgment.}  I am grateful to the referee for a careful study of
the paper and for pointing out an error in the earlier version. His or
her constructive comments lead to a significant improvement of the paper. 

This work was supported by the Hungarian Scientific
Research Fund OTKA -- T 034 512, \hbox{T 043 159}, FKFP -- 0043/2001  and the
Bolyai J\'anos Fellowship.

\section{Integrals for Hopf algebroids}
\lb{int}
\setc{0}

Hopf algebroids with bijective antipodes have been introduced in
\cite{HGD}, where several equivalent reformulations of the definition
(\cite{HGD}, Definition 4.1) have  been given. The definition
we give in this section generalizes the form in (\cite{HGD},
Proposition 4.2 {\em iii}) by allowing the antipode not to be bijective.

Integrals {\em in} Hopf algebroids also have been introduced in
\cite{HGD}. As we shall see, the definition (\cite{HGD}, Definition
5.1) applies also in our more general setting. 
In this section we introduce integrals also {\em on} Hopf algebroids.

In order for the paper to be self-contained, we recall some results on
bialgebroids from \cite{T,Lu,Sz,Sz2,KSz}. For more on bialgebroids we refer to
the papers \cite{Sch,BrzeMi,Sch1,Sz3}.

The notions of Takeuchi's $\times_R$-bialgebra \cite{T}, Lu's bialgebroid
\cite{Lu} and Xu's bialgebroid with anchor \cite{Xu} have been shown to be
equivalent in \cite{BrzeMi}. We use the definition in the following form:
\bd \label{lbgd} 
A {\em left bialgebroid} $\A_L=(A,B,s,t,\gamma,\pi)$ consists of
two algebras $A$ and $B$ over the commutative ring $k$, which are called the
total and base algebras, respectively. $A$ is a $B\stac{k}
B\op$-ring (i.e. a monoid in ${_{B\ot B\op}\M_{B\ot B\op}}$) 
via the algebra homomorphisms $s:B\to A$ and $t:B\op\to
A$, called the source and target maps, respectively. 
(This means that the ranges of $s$ and $t$ are commuting subalgebras in $A$.)
In terms of $s$ and $t$, one equips $A$ with a $B$-$B$ bimodule structure ${_B 
  A_B}$ as
$$  b\cdot a\cdot b\pri\colon = s(b) t(b\pri)a \qquad \textrm{for}\ a\in A,\
b,b\pri\in B.$$
The triple $({_B  A_B},\gamma,\pi)$ is a $B$-coring, that is a comonoid in
${_B \M_B}$. Introducing Sweedler's convention $\gamma(a)=a\di\stac{B}
a\dii$ for $a\in A$ (where implicit summation is understood), the axioms 
\bea  a\di t(b) \stac{B} a\dii &=& a\di \stac{B} a\dii s(b) \lb{cros}\\
      \gamma(1_A)&=& 1_A\stac{B} 1_A \\
      \gamma(a a\pri)&=&\gamma(a) \gamma(a\pri) \lb{gmp} \\
      \pi(1_A) &=& 1_B \\
      \pi\left(a \ \,s\!\ci\! \pi(a\pri)\right)&=&\pi(a a\pri)\lb{2.5}\\
      \pi\left(a \ \,t\!\ci\! \pi(a\pri)\right)&=&\pi(a a\pri)\lb{2.6}
\eea
are  required for all $b\in B$ and $a,a\pri\in A$.

Notice that -- although $A\stac{B} A$ is not an algebra -- axiom
(\ref{gmp}) makes sense in view of (\ref{cros}).

Homomorphisms of left bialgebroids $\A_L=(A,B,s,t,\gamma,\pi)\to 
\A_L\pri=(A\pri,B\pri,s\pri,t\pri,\gamma\pri,\pi\pri)$ are pairs of
$k$-algebra homomorphisms $(\Phi:A\to A\pri, \phi:B\to B\pri)$ satisfying
\bea s\pri\ci\phi&=&\Phi\ci s\lb{homs}\\
     t\pri\ci\phi&=&\Phi\ci t\lb{homt}\\
     \gamma \pri\ci\Phi&=&p\circ (\Phi\stac{B} \Phi)\ci \gamma\lb{homga}\\
     \pi\pri\ci\Phi&=&\phi\ci, \pi\lb{hompi}.
\eea
{
where in (\ref{homga}) $A'$ is regarded as a $B$-$B$ bimodule via $\phi$ and
$p:A'\stac B A'\to A'\stac {B'} A'$ is the canonical epimorphism.}
\ed

The bimodule ${_BA_B}$, appearing in Definition \ref{lbgd},
is defined in terms of multiplication on the left. Hence -- following the
terminology of \cite{KSz} -- we use the name {\em left} bialgebroid for this
structure. In terms of right multiplication one defines right
bialgebroids analogously. For the details we refer to \cite{KSz}.

Once the map $\gamma:A\to A\stac{B} A$ is given, we can define $\gamma\op: A\to
A\stac{B\op} A$ via $a\mapsto a\dii\ot a\di$. It is straightforward to check
that if  $\A_L=(A,B,s,t,\gamma,\pi)$ is a left bialgebroid then  $\A_{L\
cop}=(A,B\op,t,s,\gamma\op,\pi)$ is also a left bialgebroid and
$\A_L\op=(A\op,B,t,s,\gamma,\pi)$ is a right bialgebroid.

In the case of a left bialgebroid  $\A_L=(A,B,s,t,\gamma,\pi)$, the category
${_A\M}$ of left $A$-modules is
a monoidal category. As a matter of fact, any left $A$-module is a $B$-$B$
bimodule via $s$ and $t$. The monoidal product in ${_A\M}$ is defined as the
$B$-module tensor product with $A$-module structure
$$ a\cdot(m\stac{B} m\pri)\colon = a\di\cdot m\stac{B} a\dii\cdot m\pri \qquad
\textrm{ for} \ a\in A,\ m\stac{B}m\pri\in M\stac{B}M\pri . $$
Just the same way as axiom (\ref{gmp}), also this definition makes
sense in view of (\ref{cros}).
The monoidal unit is $B$ with the $A$-module structure
$$ a\cdot b\colon = \pi\left(as(b)\right)\qquad \textrm{ for} \ a\in A,\
b\in B.$$
Analogously, in the case of a right bialgebroid $\A_R$, the category $\M_A$ of
right $A$-modules is a monoidal category.

The $B$-coring structure $({_B A_B},\gamma,\pi)$, underlying a left
bialgebroid  $\A_L=(A,B,s,t,\gamma,\pi)$, gives rise to a $k$-algebra
structure 
on any of the $B$-duals of ${_B A_B}$ (\cite{BrzeWis}, 17.8). The multiplication
on the $k$-module $\asld\colon ={_B\Hom}(A,B)$, for example, is given by
\be ({_*\phi}{_*\psi})(a)={_*\psi}\left( t\!\ci\! {_*\phi}(a\dii)\
a\di\right)\qquad 
\textrm{ for}\ {_*\phi},{_*\psi}\in \asld,\ a\in A\lb{sLprod}\ee
and the unit is $\pi$.
$\asld$ is a left $A$-module and $A$ is a right $\asld$-module via
\be a\rd {_*\phi}\colon = {_*\phi}(\h a)\quad \textrm{ and}\quad 
a\ld {_*\phi}\colon = t\!\ci\! {_*\phi}(a\dii)\ a\di \lb{aslm} \ee
for ${_*\phi}\in \asld,\ a\in A$.
As it is well known \cite{Tak:sqr,KSz}, $\asld$ is also a
$B\stac{k}B\op$-ring via the inclusions
\bea {_*s}:&B\to \asld \qquad & b\mapsto \pi(\h)b \nn
     {_*t}:&B\op\to \asld \qquad & b\mapsto \pi\left(\h s(b)\right).
\nonumber\eea
Both maps ${_*s}$ and ${_*t}$ are split injections of $B$-modules with
common left 
inverse ${_*\pi}:{_*\A}\to B$, ${_*\phi}\mapsto {_*\phi}(1_A)$.
What is more, if $A$ is finitely generated and projective as a left
$B$-module, then $\asld$ has also a right bialgebroid 
structure (with source and target maps ${_*s}$ and ${_*t}$, respectively, and
counit ${_*\pi}$).

Notice that the algebra ${_*\A}$ reduces to the opposite of the usual dual
algebra if $({_B A _B},\gamma,\pi)$ is a coalgebra over a commutative ring
$B$. In the case when $A$ is a finitely generated projective left $B$-module,
also the coproduct specializes to the opposite of the usual one in the
case when $\A$ is a bialgebra. This convention is responsible
for duality to flip the notions of left- and right bialgebroids.

Applying the above formulae to the left bialgebroid $(\A_L)_{cop}$, we obtain a
$B\stac{k} B\op$-ring structure on $\atld\colon ={\Hom_B}(A,B)$. The
inclusions $B\to \atld$ and $B\op\to \atld$ will be denoted by ${s_*}$ and
${t_*}$, respectively. In particular, $\atld$ is a left $A$-module
and  $A$ is a right $\atld$-module via 
\be  a\ru {\phi_*}\colon =(\h a)\quad \textrm{ and}\quad 
a\lu{\phi_*}\colon = s\!\ci\!{\phi_*}(a\di)\ a\dii. \lb{atlm}\ee
If the module ${A}$ is finitely generated and projective as a right
$B$-module then $\A_*$ is also a right bialgebroid.

In the case of a right bialgebroid $\A_R=(A,B,s,t,\gamma,\pi)$,
application of the opposite of the multiplication formula
(\ref{sLprod}) to $(\A_R)_{cop}\op$ and to $(\A_R)\op$ 
results in $B\stac{k}B\op$-ring structures on $\asrd\colon = {\Hom_B}(A,B)$ 
and $\atrd\colon ={_B\Hom}(A,B)$, respectively. We obtain  inclusions
$s^*:B\to \asrd$, $t^*:B\op\to \asrd$, ${^*s}:B\to \atrd$ and
${^*t}:B\op\to \atrd$. 

In particular, $\asrd$ and $\atrd$ are right $A$-modules and $A$ is  
a left $\asrd$-module and a left $\atrd$-module via the formulae 
\bea
\phi^*\lu a \colon =\phi^*(a\h)\quad &\textrm{ and}&\quad 
\phi^*\ru a\colon = a\uii\ t\!\ci\! \phi^*(a\ui) \lb{asrm} \\
{^*\phi}\ld a\colon ={^*\phi}(a\h)\quad &\textrm{ and}&\quad 
{^*\phi}\rd a\colon =a\ui\ s\!\ci\! {^*\phi}(a\uii) \lb{atrm}\eea
for $\phi^*\in \asrd$, ${^*\phi}\in \atrd$ and $a\in
A$. If $A$ is finitely generated and projective as a right, or as a
left $B$-module then 
the corresponding dual is also a left bialgebroid. 

Before defining the structure that is going to be the subject of the paper,
let us stop here and introduce some
notations. Analogous notations were already used in \cite{HGD}.

When dealing with a $B\stac{k}B\op$-ring $A$,
we have to face the situation that $A$ 
carries different module structures over the base algebra $B$. In this
situation the usual notation $A\stac{B} A$ would be ambiguous. Therefore we
make the following notational convention. In terms of the maps
$s:B\to A$ and $t:B\op\to A$, we introduce four $B$-modules
\bea {_B A}:&\qquad &b\cdot a \colon = s(b)a\nn
A_B:&\qquad &a\cdot b\colon = t(b)a\nn
A^B:&\qquad&a\cdot b=as(b)\nn
{^B A}:&\qquad& b\cdot a= at(b).
\lb{amod}\eea
(Our notation can be memorized as left indeces stand for left modules and
right indeces stand for right modules. Upper indeces are used to label modules
defined in terms of right multiplication and lower indeces are used for
modules defined in terms of left multiplication.) 

In writing $B$-module tensor products, we write out explicitly the module
structures of 
the factors that are taking part in the tensor products, and do not put marks
under the symbol $\ot$. E.g. we write $A_B\ot {_B A}$. Normally
we do not denote the module structures that are not taking part in the tensor
product, this should be clear from the context. In writing elements of
tensor product modules we do not distinguish between the various module
tensor products. That is, we write both $a\ot a\pri\in A_B\ot {_B A}$ and $c\ot
c\pri\in A^B\ot {_B A}$, for example.

A left $B$-module can be considered as a right $B\op$-module, and sometimes we
want to take a module tensor product over $B\op$. In this case we use the
name of the corresponding $B$-module and the fact that the tensor product is
taken 
over $B\op$ should be clear from the order of the factors. For example, ${_B
  A}\ot A_B$ is the $B\op$-module tensor product of the right $B\op$ module
defined via multiplication by $s(b)$ on the left, and the left $B\op$-module
defined via multiplication by $t(b)$ on the left. 

In writing multiple tensor products, we use different types of letters to
denote which module
structures take part in the same tensor product. For example, the $B$-module
tensor product $A_B \ot {^B A}$ can be given a right $B$ module structure via
multiplication by $t(b)$ on the left in the second factor. The tensor product
of this right $B$-module with ${_B A}$ is denoted by
$A_B\ot {^B A_{\bf B}}\ot {_{\bf B} A}$.

\smallskip

We are ready to introduce the structure that is the main subject of the paper: 
\bd \lb{hgd}
A {\em Hopf algebroid} $\A=(\A_L,\A_R,S)$ consists of a left bialgebroid
$\A_L=(A,L,s_L,$ $t_L,\gamma_L,\pi_L)$, a right bialgebroid
$\A_R=(A,R,s_R,t_R,\gamma_R,\pi_R)$ on the {\em same} total algebra $A$, and a
$k$-module map $S:A\to A$, called the 
antipode, such that the following axioms hold true:
\bea i)&& s_L\ci \pi_L\ci t_R=t_R,\qquad t_L\ci \pi_L\ci s_R=s_R
\quad \textrm{and}\nn
&& s_R\ci \pi_R\ci t_L=t_L,\qquad t_R\ci \pi_R\ci s_L=s_L
\lb{hgdi}\\
ii)&&(\gamma_L\ot {^R A})\ci \gamma_R=(A_L\ot \gamma_R)\ci\gamma_L
\quad \textrm{as\ maps\ } A\to A_L\ot {_L A^R}\ot {^R A}\quad
\textrm{and}\nn
&&(\gamma_R\ot {_L A})\ci \gamma_L=(A^R\ot \gamma_L)\ci\gamma_R
\quad \textrm{as\ maps\ } A\to A^R\ot {^R A_L}\ot {_L A}\lb{hgdii}\\
iii)&& S\ \textrm{is\ both\ an\ }L\textrm{-}L\ \textrm{bimodule\ map\ } {^L
A_L}\to {_L A^L}\ \textrm{and\ an\ } R\textrm{-}R\ \textrm{bimodule\ map}\nn
&& {^R A_R}\to
{_R A^R} \lb{hgdiii}\\
iv)&&m_A\ci (S\ot {_L A})\ci \gamma_L=s_R\ci\pi_R\quad \textrm{and}\nn
   &&m_A\ci (A^R\ot S)\ci \gamma_R=s_L\ci \pi_L. \lb{hgdiv}
\eea
\ed
If $\A=(\A_L,\A_R,S)$ is a Hopf algebroid then so is
$\A\op_{cop}=((\A_R)\op_{cop},(\A_L)\op_{cop},S)$ and if $S$ is
bijective then also $\A_{cop}=((\A_L)_{cop},(\A_R)_{cop},S\inv)$ and
$\A\op=((\A_R)\op,(\A_L)\op,S\inv)$. 

The following modification of Sweedler's convention will turn out to
be useful. For a Hopf algebroid $\A=(\A_L,\A_R,S)$ we
use the notation $\gamma_L(a)=a\di\ot a\dii$ with lower indices, and
$\gamma_R(a)=a\ui\ot a\uii$ with upper indices for $a\in A$ in the
case of the coproducts of $\A_L$ and of
$\A_R$, respectively. The axioms (\ref{hgdii}) read in this notation
as 
\bea {a\ui}\di\ot {a\ui}\dii\ot a\uii&=&a\di\ot {a\dii}\ui\ot
{a\dii}\uii\nn
{a\di}\ui\ot {a\di}\uii\ot a\dii&=&a\ui\ot {a\uii}\di \ot {a\uii}\dii
\nonumber \eea
for $a\in A$.

Examples of Hopf algebroids (with bijective antipode) are collected in
\cite{HGD}.
\bp\lb{antip}
1) The base algebras $L$ and $R$ of the left and right bialgebroids in a
Hopf algebroid are anti-isomorphic.

2) For a Hopf algebroid $\A=(\A_L,\A_R,S)$, the pair $(S,\pi_L\ci s_R)$ is
a left bialgebroid homomorphism $(\A_R)\op_{cop}\to \A_L$ and
$(S,\pi_R\ci s_L)$ is a left bialgebroid homomorphism $\A_L\to
(\A_R)\op_{cop}$. 
\ep
\pr ${\underline{1)}}$: Both $\pi_R\ci s_L$ and $\pi_R\ci t_L$ are
anti-isomorphisms $L\to R$ with inverses $\pi_L\ci t_R$ and $\pi_L\ci s_R$,
respectively.   

${\underline{2)}}$: By part 1), the map $\pi_L\ci
s_R:R\op\to L$ is an algebra homomorphism. 
It follows from (\ref{hgdiii}), (\ref{hgdiv}) and some bialgebroid
identities that $S:A\op\to A$ is an algebra homomorphism, as for
$a,b\in A$ we have
\bea S(1_A)&=&1_A\ S(1_A)=s_L\ci\pi_L(1_A)=1_A\quad \qquad\textrm{and}\nn
S(ab)&=&S[t_L\!\ci\! \pi_L(a\dii)\ a\di \ b]\nn
&=&S[a\di\  t_L\!\ci\! \pi_L(b\dii)\ b\di]\ {a\dii}\ui S({a\dii}\uii)\nn
&=& 
{
S(a\di b\di){a\dii}\ui {b\dii}\ui S({b\dii}\uii) S({a\dii}\uii)}\nn
&=&S[{a\ui}\di {b\ui}\di]\ {a\ui}\dii {b\ui}\dii S({b}\uii) S(a\uii)\nn
&=&s_R\!\ci\! \pi_R(a\ui b\ui)\ S(b\uii)\ S(a\uii)\nn
&=&S\left[b\uii \ t_R\!\ci\! \pi_R\left(t_R\!\ci\!
  \pi_R(a\ui)\ b\ui\right)\right]S(a\uii)\nn 
&=& S(b)\ s_R\!\ci\! \pi_R(a\ui)\  S(a\uii)=S(b)S(a).
\nonumber\eea
{
In the verification of the anti-multiplicativity
of $S$, the third equality follows by axioms (\ref{cros}) and
(\ref{hgdiv}). In the fourth equality we used that, for any $x,y\in A$, there
are well defined maps $A_L \ot {}_L A^R \ot {}^R A\to A$, $a\ot b \ot
c \mapsto S(ax)b y\ui S(y\uii) S(c)$ and $a\ot b \ot c \mapsto S(x\di
a) x\dii b S(c) S(y)$, that can be composed with the equal maps 
$(\gamma_L \ot {}^R A)\circ \gamma_R = (A_L \ot \gamma_R)\circ \gamma_L: A \to
A_L \ot {}_L A^R \ot {}^R A$, cf. (\ref{hgdii}). 
The fifth equality follows by (\ref{gmp}) and (\ref{hgdiv}), 
the sixth equality is a consequence of (\ref{hgdiii}) and the right
bialgebroid analogue of (\ref{2.6}).
The last two equalities both follow by the counitality of $\gamma_R$ and
(\ref{hgdiii}), taking into account that $\gamma_R(t_R(r)b)=t_R(r)\ui b\ui
\ot t_R(r)\uii b\uii = t_R(r) b\ui \ot b\uii$, for all $r\in R$ and $b\in A$.}

Properties (\ref{homs}-\ref{homt}) follow from (\ref{hgdiii}) and
(\ref{hgdi}) as
\bea s_L\ci\pi_L\ci s_R&=&S\ci t_L\ci \pi_L\ci s_R=S\ci s_R\lb{S&L}\\
t_L\ci \pi_L\ci s_R &=& s_R=S\ci t_R. \lb{S&R}
\eea

Properties (\ref{homga}-\ref{hompi}) are checked on an element
$a\in A$ as
{
\bea \gamma_L\ci S(a) &=& 
S({a\ui}\di)\di {a\ui}\dii S(a\uii)\ot S({a\ui}\di)\dii\nn
&=& S({a\ui}\di)\di {{{a\ui}\dii}\ui}\di S(a\uii)\ot 
S({a\ui}\di)\dii  {{{a\ui}\dii}\ui}\dii S({{a\ui}\dii}\uii)\nn
&=& {\left(S({a\ui}\di){{{a\ui}\dii}\ui}\right)}\di S(a\uii)\ot
{\left(S({a\ui}\di){{{a\ui}\dii}\ui}\right)}\dii S({{a\ui}\dii}\uii)\nn
&=&(S\ot S)\ci \gamma_R\op(a) \lb{S&gamma}\\
\pi_L\ci S(a)&=& \pi_L[S(a\di)\ s_L\!\ci\! \pi_L(a\dii)]=
\pi_L[S(a\di) a\dii]=
\pi_L\ci s_R\ci\pi_R(a) \lb{S&pi}. 
\eea
In order to check the compatibility condition between $S$ and the coproducts,
note that there is a well defined map $A_L \ot {}_L
A^R \ot {}^R A\to A_L \ot {}_L A$, $a\ot b \ot c \mapsto a\ot
 bS(c)$. Composing it with the equal maps $(\gamma_L \ot {}^R A)\circ
\gamma_R = (A_L \ot \gamma_R)\circ \gamma_L: A \to A_L \ot {}_L A^R \ot {}^R
A$ (cf. (\ref{hgdii})), and using (\ref{hgdiv}), we conclude that 
\begin{equation}\label{eq:5}
{a\ui}\di \ot a {\ui}\dii S(a\uii)=a\ot 1_A, \qquad \textrm{for all
}a\in A. 
\end{equation}
Applying to both sides of (\ref{eq:5}) the well defined map $A_L\ot {}_LA \to
A_L \ot {}_L A$, $a\ot b \mapsto S(a)\di b \ot S(a)\dii$, we conclude on
the first equality of the computation in (\ref{S&gamma}). 
In the second equality we used 
(\ref{eq:5})
again. The third equality follows
by multiplicativity of $\gamma_L$, cf. (\ref{gmp}).
The last equality is derived similarly to the first one:
There is a well defined map $A_L \ot {}_L A^R \ot {}^R A \to A^R \ot {}_R A$,
$a\ot b \ot c \mapsto S(a)b \ot S(c)$. Composing it with the equal
maps $(\gamma_L \ot {}^R A)\circ \gamma_R = (A_L \ot \gamma_R)\circ \gamma_L:
A \to A_L \ot {}_L A^R \ot {}^R A$, using (\ref{hgdiv}) and the identity $S
\circ t_R =s_R$, we conclude that $S(a\di){a\dii}\ui \ot S({a\dii}\uii)=
1\ot S(a)$, for all $a\in A$. Applying to both sides of this identity the
well defined map $A^R\ot {}_R A \to A_L \ot {}_L A$, $a\ot b \mapsto a\di
\ot a\dii b$, we obtain $1_A \ot S(a)=
\left(S(a\di){a\dii}\ui\right)\di \ot \left(S(a\di){a\dii}\ui\right)\dii
S({a\dii}\uii)$, that explains the last equality in (\ref{S&gamma}).

In the first equality of (\ref{S&pi}), we used the counitality of $\gamma_L$
and (\ref{hgdiii}). The second equality follows by (\ref{2.5}). To derive the
last equality, we made use of (\ref{hgdiv}).}

The proof is completed by the observation that in passing from the Hopf
algebroid $\A$ to $\A\op_{cop}$ the roles of $(S,\pi_L\ci s_R)$ and
$(S,\pi_R\ci s_L)$ become interchanged. 
\hfill\qed

\bp\lb{schau}
The left bialgebroid $\A_L$ in a Hopf algebroid $\A$ is a
$\times_L$-Hopf algebra in the sense of \cite{Sch}. That is, the map
$$\alpha: {^L A}\ot {A_L}\to A_L\ot {_L A}\qquad a\ot b\to a\di \ot
a\dii b $$
is bijective.
\ep
\pr The inverse of $\alpha$ is given by
$$ \alpha\inv:A_L\ot {_L A}\to {^L A}\ot {A_L} \qquad a\ot b\mapsto
a\ui\ot S(a\uii)b. $$

\vspace{-.8cm}
\hfill\qed

\vspace{.2cm}

The relation between the left and the right bialgebroids in a Hopf
algebroid $\A$ implies relations between the dual algebras
$\asrd\equiv \Hom_R(A^R,R)$ and $\atld\equiv \Hom_L(A_L,L)$ and also
between $\atrd\equiv {_R \Hom}({^R A},R)$ and $\asld\equiv {_L
\Hom}({_L A},L)$:
\bl\lb{sigmachi}
For a Hopf algebroid $\A$, there exist algebra anti-isomorphisms
$\sigma:\asld\to \atrd$ and $\chi:\asrd\to \atld$ satisfying
\bea a\ld{_*\phi} &=& \sigma({_*\phi})\rd a \qquad
\textrm{and}\lb{sigma}\\
\phi^*\ru a &=& a\lu \chi(\phi^*)\lb{chi}
\eea
for all ${_*\phi}\in \asld$, $\phi^*\in \asrd$ and $a\in A$.
\el
\pr We leave it to the reader to check that the maps
\bea \sigma:&\asld\to \atrd\qquad &{_*\phi}\mapsto \pi_R(\h\ld {_*\phi})
\quad \textrm{and} \nn
\chi:&\asrd\to \atld \qquad &\phi^*\mapsto \pi_L(\phi^*\ru \h)
\nonumber \eea are algebra anti-homomorphisms satisfying
(\ref{sigma}-\ref{chi}). The inverses are given by
\bea \sigma\inv :&\atrd\to \asld \qquad& {^*\phi}\mapsto
\pi_L({^*\phi}\rd \h)\quad \textrm{and}\nn
\chi\inv:&\atld\to \asrd \qquad& {\phi_*}\mapsto \pi_R(\h\lu {\phi_*}).
\nonumber\eea

{
\vspace{-.8cm}
\hfill\qed
\bl\lb{fgp}
{
For a Hopf algebroid $\A=(\A_L,\A_R,S)$ with a bijective antipode $S$, the
following assertions hold: 

1) The module $A_L$ is finitely generated and projective if and only if the
module ${}^RA$ is finitely generated and projective. 

2) The module ${}_L A$ is finitely generated and projective if and only if the
module $A^R$ is finitely generated and projective.}
\el}
\pr 
$\underline {1)}$:
In terms of the dual bases, $\{b_i\}\subset A$ and $\{{\beta^i_*}\}\subset
\atld$ for the module $A_L$, the dual bases, $\{k_j\}\subset A$ and
$\{{^* \kappa_j}\}\subset 
\atrd$ for the module ${^R A}$, can be constructed by the requirement that
$$ \sum_j {^* \kappa_j}\ot k_j =\sum_i \pi_R\ci t_L\ci \beta^i_*\ci S\ot
S\inv(b_i)\quad 
\textrm{as\ elements\ of\ } {\atrd_R}\ot {^R A}.$$
{
The converse implication follows by applying the same reasoning to
$\A\op_{cop}$.} 

$\underline{2)}$
follows by applying part 1)
to the Hopf algebroid $\A\op$.
\hfill\qed

\smallskip

Now we turn to the study of the notion of integrals {\em in} Hopf
algebroids. For a left bialgebroid $\A_L=(A,L,s_L,t_L,\gamma_L,\pi_L)$
and a left $A$-module $M$, the {\em invariants} of $M$ with respect to $\A_L$
are the elements of
$$ {\rm{Inv}}(M)\colon =\{\ n\in M\ \vert \ a\cdot n= s_L\ci
\pi_L(a)\cdot n\quad \forall\ a\in A \ \}. $$
Clearly, the invariants of $M$ with respect to $(\A_L)_{cop}$ coincide
with its invariants with respect to $\A_L$. The invariants of a right
$A$-module $M$ with respect to a right bialgebroid $\A_R$ are defined as
the invariants of $M$ (viewed as a left $A\op$-module) with respect to
$(\A_R)\op$. 
\bd \lb{intdef}
The {\em left integrals in a left bialgebroid} $\A_L$ are the invariants of
the left regular $A$-module with respect to $\A_L$. 

The {\em right integrals in a right bialgebroid} $\A_R$ are the invariants of
the right regular $A$-module with respect to $\A_R$. 

The {\em left/right integrals in a Hopf algebroid} $\A=(\A_L,\A_R,S)$ are
the left/right integrals in $\A_L$/$\A_R$, that is, the elements of
\bea \il(\A)&=&\{\ \ell\in A \ \ \vert\ a\ell\ =s_L\!\ci\!\pi_L(a)\ \ell \quad
\forall\ a\in A\ \}\quad \textrm{and}\nn
\ir(\A)&=&\{\ \err\in A\ \vert\  \err a=\err\  s_R\!\ci\!\pi_R(a)\quad
\forall\ a\in A\ \}.
\nonumber\eea
\ed
For any Hopf algebroid $\A=(\A_L,\A_R,S)$, we have
$\il(\A)=\ir(\A\op_{cop})$ and if $S$ is bijective then also
$\il(\A)=\il(\A_{cop})=\ir(\A^{op})$. Since for $\ell\in \il(\A)$ and
$a\in A$,
$$S(\ell)a=S[t_L\!\ci\!\pi_L(a\di)\ \ell]a\dii=S(a\di\ell)a\dii=S(\ell)\ 
s_R\!\ci\!\pi_R(a), $$
we have $S\left(\il(\A)\right)\subseteq \ir(\A)$ and, similarly,
$S\left(\ir(\A)\right)\subseteq \il(\A)$. 
\bsch \lb{intsch}
The following properties of an element $\ell\in A$ are equivalent: 
\bea 1.a)&& \ell\in \il(\A)\nn
1.b)&& S(a)\ell\ui\ot \ell\uii=\ell\ui\ot a\ell\uii \quad \
\qquad\qquad\ \forall
a\in A\nn
1.c)&& a\ell\ui\ot S(\ell\uii)=\ell\ui\ot S(\ell\uii)a \qquad \qquad\ \forall
a\in A.
\nonumber\eea
The following properties of the element $\err\in A$ are also
equivalent:
\bea 2.a)&& \err\in \ir(\A)\nn
2.b)&& \err\di\ot \err\dii S(a)=\err\di a\ot \err\dii\qquad\qquad\ \forall
a\in A\nn
2.c)&& S(\err\di)\ot \err\dii a=a S(\err\di)\ot \err\dii \qquad \quad\forall
a\in A.
\nonumber\eea
\esch

\smallskip

By  comodules over a left bialgebroid
$\A_L=(A,L,s_L,t_L,\gamma_L,\pi_L)$ we mean comodules over the
$L$-coring $({_L A_L},\gamma_L,\pi_L)$, and by comodules over a right
bialgebroid $\A_R=(A,R,s_R,t_R,\gamma_R,\pi_R)$ we mean comodules over the 
$R$-coring $({^R A^R},\gamma_R,\pi_R)$. 
The pair $({_L A},\gamma_L)$ is a left comodule, and $(A_L,\gamma_L)$ is
a right comodule over the left bialgebroid $\A_L$. Since the $L$-coring
$({_L A_L},\gamma_L,\pi_L)$ possesses a grouplike element $1_A$, 
also $(L,s_L)$ is a left comodule and $(L,t_L)$ is a right comodule over
$\A_L$ (see \cite{BrzeWis}, 28.2). Similarly, $(A^R,\gamma_R)$ and $(R,s_R)$
are right comodules, and $({^R A},\gamma_R)$ and $(R,t_R)$ are left
comodules over $\A_R$.
\bd \lb{dualint}
An {\em $s$-integral on a left bialgebroid}
$\A_L=(A,L,s_L,t_L,\gamma_L,\pi_L)$ is a left $\A_L$-comodule map
${_*\rho}:({_L A},\gamma_L)\to (L,s_L)$. That is, an element of
$$ \ir(\asld)\colon = \{\ {_* \rho}\in \asld\ \vert \ (A_L\ot
{_*\rho})\ci \gamma_L=s_L\ci {_*\rho} \ \}. $$
A {\em $t$-integral on $\A_L$} is a right $\A_L$-comodule map
$(A_L,\gamma_L)\to (L,t_L)$. That is, an element of 
$$\ir(\atld)\colon =\{\ \rho_*\in \atld\ \vert\ (\rho_*\ot {_L A})\ci
\gamma_L=t_L\ci {\rho_*}\ \}.$$
An {\em $s$-integral on a right bialgebroid}
$\A_R=(A,R,s_R,t_R,\gamma_R,\pi_R)$ is a right $\A_R$-comodule map
$(A^R,\gamma_R)\to (R,s_R)$. That is, an element of 
$$\il(\asrd)\colon =\{\ \lambda^*\in \asrd\ \vert\ (\lambda^*\ot {^R A})\ci
\gamma_R=s_R\ci \lambda^*\ \}.$$
A {\em $t$-integral on $\A_R$} is a left $\A_R$-comodule map $({^R
A},\gamma_R)\to (R,t_R)$. That is, an element of 
$$ \il(\atrd)\colon =\{\ {^*\lambda}\in \atrd\ \vert\ (A^R\ot {^*\lambda})\ci
\gamma_R=t_R\ci {^*\lambda}\ \}.$$
The {\em right/left $s$- and $t$-integrals on a Hopf algebroid}
$\A=(\A_L,\A_R,S)$ are the $s$- and $t$-integrals on
$\A_L$/$\A_R$. 
\ed
The integrals on a {\em left/right} bialgebroid are checked to be
invariants of the appropriate {\em right/left} regular module -- justifying
our usage of the terms {\em `right'} and {\em `left'} integrals for
them (cf. the remark in Section \ref{int} about our using the opposite -
co-opposite of the convention
usual in the case of bialgebras, when defining the dual
bialgebroids $\asld$ and $\asrd$). As a matter of fact, for example, if
${_*\rho}\in  \ir(\asld)$ then 
\be [{_*\rho}\ {_*\phi}](a)= {_*\phi}(a\lu {_*\rho})=
{_*\phi}(s_L\ci{_*\rho}(a))={_*\rho}(a)\
    {_*\phi}(1_A)=[{_*\rho}\ \ {_*s}\!\ci\!{_*\pi}({_*\phi})](a)
\lb{inv} \ee
for all ${_*\phi}\in \asld$ and $a\in A$.
If the module ${_L A}$ is finitely generated and projective (hence
$\asld$ is a right bialgebroid) then also the converse is true, so in this
case the $s$-integrals on $\A_L$ are the same as the right integrals
in $\asld$. Similar statements hold true on the elements of
$\ir(\atld)$, $\il(\asrd)$ and $\il(\atrd)$. 

The reader should be warned that integrals on Hopf algebras $H$
over commutative rings $k$ are defined in the literature sometimes as
comodule maps $H\to k$ -- similarly to our Definition \ref{dualint}
--, sometimes by the analogue of the weaker invariant condition
(\ref{inv}). 

For any Hopf algebroid $\A$, we have
$\ir(\asld)=\il((\A\op_{cop})^*)$ and
$\ir(\atld)=\il({^*(\A\op_{cop})})$. If the antipode is bijective then also
$\ir(\asld)=\ir((\A_{cop})_*)=\il({^*(\A\op}))$. 
\bsch \lb{dualsch}
Let $\A=(\A_L,\A_R,S)$ be a Hopf algebroid. The following
properties of an element ${_*\rho}\in \asld$ are equivalent:
\bea 1.a)&& {_*\rho}\in \ir(\asld)\nn
1.b)&& \pi_R\ci s_L\ci {_* \rho}\in \il(\atrd)\nn
1.c)&& s_L\!\ci\! {_*\rho}\left(a S(b\di)\right) \ b\dii= 
t_L\!\ci\! {_*\rho}\left(a\dii S(b)\right) \ a\di \qquad \forall a,b\in A.
\nonumber\eea
The following properties of an element $\rho_*\in \atld$ are
equivalent:
\bea 2.a)&& \rho_*\in \ir(\atld)\nn
2.b)&& \pi_R\ci t_L\ci \rho_*\in \il(\asrd)\nn
2.c)&& t_L\!\ci\! \rho_*(ab\ui) \ S(b\uii)=
s_L\!\ci\! \rho_*(a\di b)\ a\dii \qquad
 \qquad\ \ \forall a,b\in A.
\nonumber\eea
The following properties of an element $\lambda^*\in \asrd$ are
equivalent: 
\bea 3.a)&& \lambda^*\in \il(\asrd)\nn
3.b)&& \pi_L\ci s_R\ci \lambda^*\in \ir(\atld)\nn
3.c)&& a\ui \ s_R\!\ci\! \lambda^*\!\left( S(a\uii)b\right)=
b\uii \ t_R\!\ci\! \lambda^*\!\left( S(a) b\ui\right) \qquad \forall  a,b\in A.
\nonumber\eea
The following properties of an element ${^*\lambda}\in \atrd$ are
equivalent:
\bea 4.a)&& {^*\lambda}\in \il(\atrd)\nn
4.b)&& \pi_L\ci t_R\ci {^*\lambda}\in \ir(\asld)\nn
4.c)&&S(a\di)\ t_R\!\ci\! {^*\lambda} (a\dii b)= 
b\ui \ s_R\!\ci {^* \lambda} (a b\uii) \qquad  \qquad \ \forall a,b\in A.
\nonumber\eea
In particular, for ${_*\rho}\in \ir(\asld)$ the element ${_*\rho} \ci S$
belongs to $\ir(\atld)$ and for $\lambda^*\in \il(\asrd)$ the element
$\lambda^*\ci S$ belongs to $\il(\atrd)$.
\esch

\section{Maschke type theorems}
\lb{maschke}
\setc{0}

The most classical version of Maschke's theorem \cite{M} considers
group algebras over fields. It states that the group algebra of a
finite group $G$ over a field $F$ is semi-simple if and only if the
characteristic of $F$ does not divide the order of $G$. This result
has been generalized to finite dimensional Hopf algebras $H$ over 
fields $F$ by Sweedler \cite{Swe} proving that $H$ is a separable
$F$-algebra if and only if it is semi-simple and 
if and only if there exists a normalized left integral in
$H$. The proof goes as follows. It is a classical result that a
separable algebra over a field is semi-simple. If $H$ is semi-simple
then, in particular, the $H$-module on $F$, defined in terms of the
counit, is projective. This means that the counit, as an $H$-module map
$H\to F$, splits. Its right inverse maps the unit of $F$ into a
normalized integral. Finally, in terms of a normalized integral one can
construct an $H$-bilinear right inverse for the multiplication map
$H\stac{F} H\to H$. 

The only difficulty in the generalization of Maschke's theorem to Hopf
algebras over commutative rings comes from the fact that in the case of an
algebra $A$ over a commutative base ring 
$k$, separability does not imply the semi-simplicity of $A$ in the
sense \cite{Pierce} that every (left or right) $A$-module was projective. It
implies  \cite{Hattori,HirSug}, however, that every $A$-module is
$(A,k)$-projective, i.e. that every epimorphism of
$A$-modules which is $k$-split, is also $A$-split. 
In order to avoid confusion, we will say that the $k$-algebra $A$ is
{\em semi-simple} \cite{Pierce} if it is an Artinian semi-simple ring i.e. if
any $A$-module is 
projective. By the terminology of \cite{Hattori} we call $A$ a (left or
right) {\em semi-simple extension} of $k$ if any (left or right) $A$-module
is $(A,k)$-projective.

Since the counit of a Hopf algebra $H$ over a commutative ring $k$ is a split
epimorphism of $k$-modules, Maschke's theorem generalizes to this case in
the following form \cite{CaenMi,Lomp}. The extension $k\to H$ is
separable if and only if it is (left and right) semi-simple and if and
only if there exist normalized (left and right) integrals in $H$. 

In this section we investigate the properties of the total algebra of
a Hopf algebroid as an extension of the base algebra, that are equivalent to
the 
existence of normalized integrals {\em in} the Hopf algebroid. Dually,
we investigate also the properties of the coring over the base algebra
underlying a Hopf algebroid, that are equivalent to the existence of 
normalized integrals {\em on} the Hopf algebroid (in any of the four possible
senses). 

A Maschke type theorem on certain Hopf algebroids can be obtained also by 
application of (\cite{Sz4}, Theorem 4.2). Notice, however, that the Hopf
algebroids occurring this way are only the Frobenius Hopf algebroids
(discussed in Section \ref{frob} below), that is the Hopf algebroids
possessing non-degenerate integrals (which are called Frobenius integrals in
\cite{Sz4}). 

\bigskip

Following Theorem \ref{maschkethm} generalizes results from
(\cite{CaenMi}, Proposition 4.7) and (\cite{Lomp}, Theorem 3.3). 
\bt \lb{maschkethm} (Maschke Theorem for Hopf algebroids.)
The following assertions on a Hopf algebroid $\A=(\A_L,\A_R,S)$ are
equivalent: 

1.a) The extension $s_R:R\to A$ is separable. That is, the multiplication map
   $A^R\ot {_R A}\to A$ splits as an $A$-$A$ bimodule epimorphism.

1.b) The extension $t_R:R\op\to A$ is separable. That is, the multiplication map
   ${^R A}\ot {A_R}\to A$ splits as an $A$-$A$ bimodule epimorphism.

1.c) The extension $s_L:L\to A$ is separable. That is, the multiplication map
   $A^L\ot {_L A}\to A$ splits as an $A$-$A$ bimodule epimorphism.

1.d) The extension $t_L:L\op\to A$ is separable. That is, the multiplication map
   ${^L A}\ot {A_L}\to A$ splits as an $A$-$A$ bimodule epimorphism.

2.a) The extension $s_R:R\to A$ is right semi-simple. That is, any right
   $A$-module is $(A,R)$-projective.

2.b) The extension $t_R:R\op\to A$ is right semi-simple. That is, any right
   $A$-module is $(A,R\op)$-projective.

2.c) The extension $s_L:L\to A$ is left semi-simple. That is, any left
   $A$-module is $(A,L)$-projective.

2.d) The extension $t_L:L\op\to A$ is left semi-simple. That is, any left
   $A$-module is $(A,L\op)$-projective.

3.a) There exists a normalized left integral in $ \A$.
That is, an element $\ell\in\il(\A)$ such that\break
$\pi_L(\ell)=1_L$.

3.b) There exists a normalized right integral in $ \A$.
That is, an element $\err\in\ir(\A)$ such that $\pi_R(\err)=1_R$.

4.a) The\ epimorphism $\pi_R:A\to R$ splits as a right $A$-module map.

4.b) The\ epimorphism $\pi_L:A\to L$ splits as a left $A$-module map.
\et
\pr $\underline {1.a)\Rightarrow 2.a)}$, $\underline {1.b)\Rightarrow
2.b)}$, $\underline {1.c)\Rightarrow 2.c)}$ and $\underline
{1.d)\Rightarrow 2.d)}$:
It is proven in (\cite{HirSug}, Proposition 2.6) that a separable
extension is both left- and right semi-simple.

$\underline{2.a)\Rightarrow 4.a)}$ ( $\underline{2.b)\Rightarrow 4.a)}$ ):
The epimorphism $\pi_R$ is split as a right ( left ) $R$-module map by
$s_R$ ( by $t_R$ ), hence it is split as a right $A$-module map.

$\underline{4.a)\Rightarrow 3.b)}$: Let $\nu: R\to A$ be the right inverse
of $\pi_R$ in ${\M_A}$. Then $\err\colon =\nu(1_R)$ is a normalized
right integral in $\A$.

$\underline{3.a)\Leftrightarrow 3.b)}$: By part 2) of Proposition \ref{antip}
the antipode takes a normalized 
left/right integral to a normalized right/left integral.

$\underline{3.a)\Rightarrow 1.a)}$ and $\underline{3.b)\Rightarrow
1.b)}$: 
If $\ell$ is a normalized left integral in $\A$ then, by Scholium
\ref{intsch}, the required right inverse of the multiplication map $A^R\ot{_R
  A}\to A$ is given by the $A$-$A$ bimodule map $a\mapsto a\ell\ui\ot
S(\ell\uii)\equiv \ell\ui\ot S(\ell\uii)a$. Similarly, if $\err$ is a
normalized right integral in $\A$ then the right inverse of the multiplication
map ${^R A}\ot A_R\to A$ is given by $a\mapsto a S(\err\di)\ot \err\dii\equiv
S(\err\di)\ot \err\dii a$.

The proof is completed by applying the above arguments to the Hopf
algebroid $\A\op_{cop}$. \hfill\qed
\smallskip

Let us make a comment on the semi-simplicity of the algebra $A$
(cf. \cite{HirSug}, Proposition 1.3). If $R$ is a 
semi-simple algebra and the equivalent conditions of Theorem
\ref{maschkethm} hold true, then $A$ -- being a semi-simple extension of
a semi-simple algebra -- is a semi-simple algebra. On the other hand,  
notice that  condition $4.a)$ in Theorem \ref{maschkethm}
is equivalent to the projectivity of the right $A$-module 
$R$. Hence if $A$ is a semi-simple $k$-algebra then the equivalent
conditions of the theorem hold true. It is not true, however, that the
semi-simplicity of the total algebra implied the semi-simplicity of the base
algebra (which was shown by Lomp to be the case in Hopf algebras
\cite{Lomp}). A counterexample can be constructed as follows: If $B$ is a
Frobenius algebra over a commutative ring $k$ then $A\colon = \End_k(B)$ has a
Hopf algebroid structure over the base $B$ \cite{Frob}. If $B$ is a
Frobenius algebra over a field -- which can be non-semi-simple! --
then $A$ is a Hopf algebroid with semi-simple total algebra.

\smallskip

Following Theorem \ref{dualmthm} can be considered as a dual of
Theorem \ref{maschkethm} in the sense that it speaks about corings over the
base algebras instead of algebra extensions. It is
important to emphasize, however, that the two theorems are
independent results. Even in the case of Hopf algebroids such that all
module structures (\ref{amod}) are finitely generated and projective,
the duals are not known to be Hopf algebroids. 

Recall that the dual notion of that of a relative projective module is
the relative injective comodule. Namely, a comodule $M$ for an $R$-coring $A$
is called {\em $(A,R)$-injective} (\cite{BrzeWis}, 18.18) if any monomorphism of
$A$-comodules from $M$, which splits as an $R$-module map, splits also as an
$A$-comodule map. 

\bt \lb{dualmthm} (Dual Maschke Theorem for Hopf algebroids.)
The following assertions on a Hopf algebroid $\A=(\A_L,\A_R,S)$ are
equivalent:

1.a) The $R$-coring $({^R A^R},\gamma_R,\pi_R)$ is coseparable. That is, 
the comultiplication $\gamma_R:A\to A^R\ot {^R A}$ splits as an
$\A_R$-$\A_R$ bicomodule monomorphism. 

1.b) The $L$-coring $({_LA_L},\gamma_L,\pi_L)$ is coseparable.  That is, the
comultiplication $\gamma_L:A\to A_L\ot {_L A}$ splits as an
$\A_L$-$\A_L$ bicomodule monomorphism.  

2.a) Any right $\A_R$-comodule is $(\A_R,R)$-injective.

2.b) Any left $\A_R$-comodule is $(\A_R,R)$-injective.

2.c) Any left $\A_L$-comodule is $(\A_L,L)$-injective.

2.d) Any right $\A_L$-comodule is $(\A_L,L)$-injective.

3.a) There exists a normalized left $s$-integral on $\A$.
That is, an element $\lambda^*\in\il(\asrd)$ such that
$\lambda^*(1_A)=1_R$.

3.b) There exists a normalized left $t$-integral on $\A$. That is, an
element ${^*\lambda}\in\il(\atrd)$ such that
${^*\lambda}(1_A)=1_R$.

3.c) There exists a normalized right $s$-integral on $\A$.
That is, an element ${_*\rho}\in\ir(\asld)$ such\ that
${_*\rho}(1_A)=1_L$.

3.d) There exists a normalized right $t$-integral on $\A$. That is, an
element $\rho_*\in\ir(\atld)$ such that $\rho_*(1_A)=1_L$.

4.a) The monomorphism $s_R:R\to A$ splits as a right ${\A_R}$-comodule map.
 
4.b) The monomorphism $t_R:R\to A$ splits as a left ${\A_R}$-comodule map.

4.c) The monomorphism $s_L:L\to A$ splits as a left ${\A_L}$-comodule map.

4.d) The monomorphism $t_L:L\to A$ splits as a right ${\A_L}$-comodule map.
\et
\pr
$\underline{1.a)\Rightarrow 2.a), 2.b)}$ is proven in (\cite{BrzeWis}, 26.1).

$\underline{2.a)\Rightarrow 4.a)}$ ( $\underline{2.b)\Rightarrow 4.b)}$ ): 
The monomorphism $s_R$ ( $t_R$ ) is split as a right ( left ) $R$-module
map by $\pi_R$ hence it is split as a right ( left ) $\A_R$-comodule map.

$\underline{4.a)\Rightarrow 3.a)}$ and $\underline{4.b)\Rightarrow 3.b)}$: The
left 
inverse $\lambda^*$ of $s_R$ in the category of right ${\A_R}$-comodules  is a
normalized $s$-integral on $\A_R$ by very definition. Similarly, the left
inverse ${^*\lambda}$ of $t_R$ in the category of left ${\A_R}$-comodules is a
normalized $t$-integral on $\A_R$. 

$\underline{3.a)\Rightarrow 3.b)}$: If $\lambda^*$ is a normalized
$s$-integral on $\A_R$ then $\lambda^*\ci S$ is a normalized $t$-integral on
$\A_R$ by Scholium \ref{dualsch}.

$\underline{3.b)\Rightarrow 1.a)}$: 
In terms of the normalized $t$-integral ${^*\lambda}$ on $\A_R$ the
required right inverse of the coproduct $\gamma_R$ is constructed as
the map
$$A^R\ot {^R A}\to A,\qquad a\ot b\mapsto 
t_R\!\ci\! {^*\lambda}\left(aS(b\di)\right)\ b\dii. $$
It is checked to be an $\A_R$-$\A_R$ bicomodule map using that by
Scholium \ref{dualsch}, 4.b) and 1.c) we have $t_R\!\ci\!
{^*\lambda}\left(aS(b\di)\right)\ b\dii= a\ui\
s_R\!\ci\!\pi_R[\ t_R\!\ci\!{^* \lambda}\left(a\uii S(b\di)\right)\
b\dii\ ]$ for all $a,b$ in $A$.

$\underline{3.a)\Leftrightarrow 3.d)}$ follows from Scholium
\ref{dualsch}, 2.b).

\noindent
The remaining equivalences are proven by applying the above arguments to the
Hopf algebroid $\A\op_{cop}$. 
\hfill\qed
\smallskip

The proofs of Theorem \ref{maschkethm} and \ref{dualmthm} can be unified
if one formulates them as equivalent statements on the forgetful functors from
the category of $A$-modules, and from the category of $\A_L$ or
$\A_R$-comodules, respectively, to the category of $L$- or $R$-modules -- as
it is done in the case of Hopf algebras over commutative rings 
in \cite{CaenMi}. We belive (together with the referee), however, that the
above formulation in terms of algebra extensions and corings, respectively, is
more appealing.

\section{Frobenius Hopf algebroids and non-degenerate integrals}
\lb{frob}
\setc{0}

A left or right integral  $\ell$ in a Hopf algebra
$(H,\Delta,\epsilon,S)$ over a 
commutative ring $k$ is called non-degenerate \cite{LaSwe} if the maps 
\bea \Hom_k(H,k)\to H&\qquad& \phi\mapsto (\phi\ot H)\ci \Delta(\ell)
\quad\textrm{and}\nn
\Hom_k(H,k)\to H&\qquad& \phi\mapsto (H\ot \phi)\ci \Delta(\ell)
\nonumber\eea
are bijective. 

The notion of non-degenerate integrals is made relevant by the Larson-Sweedler
Theorem \cite{LaSwe} stating that a free and finite dimensional
bialgebra over a principal ideal domain is a Hopf algebra if and only if there
exists a non-degenerate left integral in $H$.

The Larson-Sweedler Theorem has been extended by Pareigis \cite{Par} to Hopf
algebras over commutative rings with trivial Picard group.
He proved also that a bialgebra over an arbitrary commutative ring
$k$, which is 
a Frobenius $k$-algebra, is a Hopf algebra if and only if
there exists a Frobenius functional $\psi:H\to k$ satisfying in addition
$$ (H\ot \psi)\ci \Delta =1_H\psi(\h). $$
As a matter of fact, based on the results of \cite{Par},
the following variant of (\cite{CaMiZhu}, 3.2 Theorem 31) can be proven:
\bt \lb{parthm}
The following properties of a Hopf algebra $(H,\Delta,\epsilon,S)$ over a
commutative ring $k$ are equivalent:

1) $H$ is a Frobenius $k$-algebra.

2) There exists a non-degenerate left integral in $H$.

3) There exists a non-degenerate right integral in $H$.

4) There exists a non-degenerate left integral on $H$. That is, a
Frobenius functional $\psi: H\to k$ satisfying 
$ (H\ot \psi)\ci \Delta=1_H\psi(\h)$.
  
5) There exists a non-degenerate right integral on $H$. That is, a
Frobenius functional $\psi: H\to k$ satisfying 
$ (\psi\ot H)\ci \Delta=1_H\psi(\h)$.
\et
The main subject of the present section is a generalization of 
Theorem \ref{parthm} to Hopf algebroids.

\smallskip

The most important tool in the proof of Theorem \ref{parthm} is the Fundamental
Theorem for Hopf modules \cite{LaSwe}. A very general form of it has
been proven by Brzezi\'nski (\cite{Brze}, Theorem 5.6, see also
\cite{BrzeWis}, 28.19) in
the framework of corings. It can be applied in our setting as follows.

Hopf modules over bialgebroids are examples of Doi-Koppinen modules over
algebras, studied in \cite{BrzeCaenMil}. A left-left Hopf module over a left
bialgebroid $\A_L=(A,L,s_L,t_L,\gamma_L,\pi_L)$ is a left comodule for the
comonoid $(A,\gamma_L,\pi_L)$ in the category of left
$A$-modules. That is, a pair $(M,\tau)$ where $M$
is a left $A$-module, hence a left $L$-module ${_L M}$ via $s_L$. The pair
$({_L M},\tau)$ is a left $\A_L$-comodule such that $\tau:M\to A_L\ot {_L M}$
is a left $A$-module map to the module
$$ a\cdot (b\ot m)\colon = a\di b\ot a\dii\cdot m\qquad \textrm{for}\ a\in
A,\ b\ot m\in  A_L\ot {_L M}.$$
The right-right Hopf modules over a right bialgebroid $\A_R$ are the left-left
Hopf modules over $(\A_R)\op_{cop}$.

It follows from (\cite{BrzeCaenMil}, Proposition 4.1) that the left-left Hopf
modules over $\A_L$ are the left comodules over the $A$-coring
\be {\cal W}\colon =(A_L\ot {_L A}, \gamma_L\ot{_L A}, \pi_L\ot {_L A}),
\lb{w}\ee 
where the $A$-$A$ bimodule structure is given by
$$ a\cdot (b\ot c)\cdot d\colon = a\di b\ot a\dii cd\qquad \textrm{for}\
a,d\in A,\ b\ot c\in A_L\ot {_L A}.$$
The coring (\ref{w}) was studied in \cite{int}. It was shown to possess
a group-like element $1_A\ot 1_A\in A_L\ot {_L A}$ and corresponding
coinvariant 
subalgebra $t_L(L)$ in $A$. The coring (\ref{w}) is Galois (w.r.t. the
group-like element  $1_A\ot 1_A$) if and only if $\A_L$ is a $\times_L$-Hopf
algebra in the sense of \cite{Sch}. Since in a Hopf algebroid
$\A=(\A_L,\A_R,S)$ the left bialgebroid $\A_L$ is a $\times_L$-Hopf algebra by
Proposition \ref{schau}, the $A$-coring (\ref{w}) is Galois in this
case. Denote the category of left-left Hopf modules over $\A_L$
(i.e. of left comodules over the coring (\ref{w})) by ${^{\cal W}\M}$. The
application of (\cite{Brze}, Theorem 5.6) results that if
$\A=(\A_L,\A_R,S)$ is a Hopf algebroid, such that the module ${^L A}$ 
is faithfully flat, then the functor 
\be G: {^{\cal W}\M}\to \M_L\qquad (M,\tau)\mapsto \textrm{Coinv}(M)_L\colon =\{\
  m\in M\ \vert\ \tau(m)=1_A\ot m\in A_L\ot {_L M}\ \}\lb{G}\ee
(where the right $L$-module structure on $\textrm{Coinv}(M)$ is given via
$t_L$) and the induction functor
\be F: \M_L\to {^{\cal W}\M}\qquad N_L\mapsto ({^L A}\ot N_L,\gamma_L\ot N_L)
\lb{F}\ee
(where the left $A$-module structure on ${^L A}\ot N_L$ is given by left
multiplication in the first factor) are inverse equivalences.

In the case of Hopf algebras $H$ over commutative rings $k$, these
arguments lead to the Fundamental Theorem only for faithfully flat Hopf
algebras. The proof of the Fundamental Theorem in \cite{LaSwe},
however, does not rely on any assumption on the $k$-module structure
of $H$. 
Since the Hopf algebroid structure is more restrictive than the 
$\times_L$-Hopf algebra structure, one hopes to prove the
Fundamental Theorem for Hopf algebroids also under milder assumptions
-- using the whole strength of the Hopf algebroid structure.  

{
In the following theorem, Sweedler's index notation
$\tau(m)=m_{\langle-1\rangle}\ot {m_{\langle0\rangle}}$ (with implicit
summation) is used, for the left coaction $\tau:M\to {\cal W}\ot_A M \cong 
A_L \ot {}_L M$ of the constituent left $L$-bialgebroid ${\mathcal A}_L$ in a
Hopf algebroid ${\mathcal A}$, on a left ${\mathcal A}_L$-comodule $M$ and
$m\in M$.} 

\bt\lb{fundi} (Fundamental Theorem for Hopf algebroids.)
Let $\A=(\A_L,\A_R,S)$ be a Hopf algebroid and ${\cal W}$ be the $A$-coring
(\ref{w}). 
{
Assume that the kernel of the maps
\begin{equation}\label{eq:kernel}
M \to A_L \ot {}_L M,\qquad 
m\mapsto ({m_{\langle-1\rangle}} \ot {m_{\langle 0 \rangle}}) - (1_A \ot m)
\end{equation}
is preserved by the
functor ${}^L A \ot - : \M_L \to \M_L$, for any $M\in {^{\cal W}\M}$
(e.g. ${}^LA$ is a flat module).}
Then the functors $G: {^{\cal W} \M}\to \M_L$ in (\ref{G}) and
$F:\M_L\to {^{\cal W}\M}$ in (\ref{F}) are inverse equivalences.\footnote{
{
In the arXiv version of \cite{Bohm:rev}, a more restrictive notion of a
comodule of a Hopf algebroid is studied, cf. \cite[arXiv version, Definition
  2.19]{Bohm:rev}. The total algebra $A$ of any Hopf algebroid ${\cal A}$ can
be regarded as a monoid in the monoidal category of ${\cal A}$-comodules in
this more 
restrictive sense. In this setting, the category of $A$-modules in the
category of ${\cal A}$-comodules, and the category of modules for the base
algebra $L$ of ${\cal A}$, were proven to be equivalent, without any further
(equalizer preserving) assumption, see \cite[Theorem 3.26 and Remark
  3.27]{Bohm:rev}. That is, in the arXiv version of \cite{Bohm:rev}, a weaker
statement is proven under weaker assumptions.
}}
\et
\pr We construct natural isomorphisms $\alpha:F\ci G\to {^{\cal W}\M}$ and
$\beta:G\ci F\to \M_L$. The map
$$\alpha_M:{^L A}\ot \textrm{Coinv}(M)_L\to M\qquad a\ot m\mapsto a\cdot m $$
is a left ${\cal W}$-comodule map and natural in $M$. The isomorphism property
is proven by constructing the inverse
$$\alpha_M\inv:M\to {^L A}\ot \textrm{Coinv}(M)_L\qquad m\mapsto
{m_{\langle-1\rangle}}\ui \ot S({m_{\langle-1\rangle}}\uii)\cdot
m_{\langle0\rangle}. $$   
It requires some work to check that $\alpha_M\inv(m)$ belongs to ${^L A}\ot
\textrm{Coinv}(M)_L$. 
{
By the assumption that the kernel of (\ref{eq:kernel})
is preserved
by the functor ${}^L A \ot - : \M_L \to \M_L$, we need to show only that 
\begin{equation}\label{eq:to_show}
{m_{\langle-1\rangle}}\ui \ot 
\left(S({m_{\langle-1\rangle}}\uii)\cdot
m_{\langle0\rangle}\right)_{\langle-1\rangle} \ot 
\left(S({m_{\langle-1\rangle}}\uii)\cdot
m_{\langle0\rangle}\right)_{\langle 0 \rangle} =
{m_{\langle-1\rangle}}\ui \ot 1_A \ot 
S({m_{\langle-1\rangle}}\uii)\cdot m_{\langle0\rangle},
\end{equation}
as elements of ${}^LA \ot A_{\bf L}\ot {}_{\bf L}M_L$,
for all $m\in M$.
Compose the well defined map
$$
A^R \ot {}^R A_L \ot {}_L A \to A^R \ot {}_R A, \qquad 
a\ot b \ot c \mapsto a \ot S(b)c
$$
with the equal maps $(\gamma_R \ot {}_LA)\circ \gamma_L = (A^R \ot
\gamma_L)\circ \gamma_R: A \to A^R \ot {}^R A_L \ot {}_L A$
(cf. (\ref{hgdii})) in order to conclude that, for any $a\in A$,
\begin{equation}\label{eq:key}
{a\di}\ui\ot S({a\di}\uii)a\dii
= a\ui \ot S({a\uii}\di) {a\uii}\dii
=a\ui \ot s_R \circ \pi_R(a\uii)
=a\ot 1_A.
\end{equation}
In (\ref{eq:key}), in the second equality (\ref{hgdiv}) was used, and the last
equality follows by the counitality of $\gamma_R$. Using the left
$A$-linearity of the coaction $\tau:M\to {\cal W}\ot_A M \cong 
A_L \ot {}_L M$, anti-comultiplicativity of the antipode (cf. Proposition
\ref{antip}(2)), coassociativity of $\tau$ and $\gamma_R$
and finally (\ref{eq:key}), the left hand side of (\ref{eq:to_show}) is
computed to be equal to 
\begin{eqnarray*}
{m_{\langle-2\rangle}}\ui &\ot& 
S({m_{\langle-2\rangle}}\uii)\di {m_{\langle-1\rangle}} \ot
S({m_{\langle-2\rangle}}\uii)\dii \cdot {m_{\langle 0 \rangle}}\\
&=& {m_{\langle-2\rangle}}\ui \ot 
S({{m_{\langle-2\rangle}}\uii}\uii){m_{\langle-1\rangle}} \ot
S({{m_{\langle-2\rangle}}\uii}\ui)\cdot {m_{\langle 0 \rangle}}\\
&=& {{{m_{\langle-1\rangle}}\di}\ui}\ui \ot 
S({{m_{\langle-1\rangle}}\di}\uii){{m_{\langle-1\rangle}}}\dii \ot 
S({{{m_{\langle-1\rangle}}\di}\ui}\uii)\cdot {m_{\langle 0 \rangle}}\\
&=&{m_{\langle-1\rangle}}\ui \ot 
1_A \ot
S({m_{\langle-1\rangle}}\uii)\cdot {m_{\langle 0 \rangle}}.
\end{eqnarray*}
Thus}
it follows that $\alpha_M\inv(m)$ belongs to ${^L A}\ot
\textrm{Coinv}(M)_L$ for all $m\in M$, as stated.
{
By (\ref{hgdiv}) and the counitality of $\tau$, $\alpha_M\circ
\alpha_M^{-1}(m)=m$, for all $m\in M$. It follows by (\ref{eq:key}) that also
$\alpha_M^{-1}\circ \alpha_M (a\ot m) = a\ot m$, for all $a\ot m\in {}^L A \ot
\textrm{Coinv}(M)_L$.} 

The coinvariants of the left ${\cal W}$-comodule ${^L A} \ot N_L$ are the
elements of
$$ \textrm{Coinv}({^L A} \ot N_L)=\{\ \sum_i a_i\ot n_i\in {^L A} \ot N_L\
\vert\ \sum_i a_i\ot n_i= \sum_i s_R\ci \pi_R(a_i)\ot n_i\ \}, $$
hence the map 
$$\beta_N:\textrm{Coinv}({^L A} \ot N_L)\to N\qquad \sum_i a_i\ot n_i \mapsto
\sum_i n_i\cdot \pi_L\ci S(a_i)\equiv \sum_i n_i\cdot \pi_L(a_i)$$
is a right $L$-module map and natural in $N$. It is an isomorphism with
inverse
$$\beta_N\inv:N\to \textrm{Coinv}({^L A} \ot N_L)\qquad n\mapsto 1_A\ot
n.$$

\hfill\qed 
\smallskip

An analogous result for right-right Hopf modules over $\A_R$ can be obtained by 
applying Theorem \ref{fundi} to the Hopf algebroid $\A\op_{cop}$.
\bp Let $\A=(\A_L,\A_R,S)$ be a Hopf algebroid and $(M,\tau)$ be a
left-left Hopf module over $\A_L$. Then $\textrm{Coinv}(M)$ is a
$k$-direct summand of $M$.
\ep
\pr 
{
In light of (\ref{eq:key})},
the canonical inclusion $\textrm{Coinv}(M)\to M$ is split by the
$k$-module map
\be E_M:M\to \textrm{Coinv}(M)\qquad m\mapsto S(m_{\langle-1\rangle})\cdot
m_{\langle0\rangle}. \lb{Em}\ee
\hfill\qed
\smallskip

{
As the next step towards our goal, let us assume that $\A=(\A_L,\A_R,S)$ is a
Hopf algebroid such that the module 
$A_L$ is finitely generated and projective. Under this  assumption we equip
$\asrd$ with the structure of a left-left Hopf module over $\A_L$.
Similarly, in the case when the module ${}^RA$ is finitely generated and
projective, we equip $\asrd$ with the structure of a right-right Hopf module
over $\A_R$.}  

Let $\{b_i\}\subset A$ and $\{\beta^i_*\}\subset\atld$ be dual bases for the
module $A_L$. A left $\A_L$-comodule structure on $\asrd$ can be introduced
via the $L$-action
$$ {_L \asrd}:\quad l\cdot \phi^*\colon =\phi^*\lu S\ci s_L(l)\qquad
\textrm{for}\ l\in L,\ \phi^*\in \asrd$$
and the left coaction
\be \tau_L:\asrd\to A_L\ot {_L \asrd}\qquad \phi^*\mapsto \sum_i b_i\ot
\chi\inv(\beta^i_*)\phi^* \lb{taul}.\ee
Similarly, 
{
let $\{k_j\}\subset A$ and $\{{}^*\kappa^j\}\subset\atrd$ be dual bases for
the module ${}^RA$.}
A right $\A_R$-comodule structure on $\asrd$ can be introduced
by the right $R$-action
$$ \asrd_R:\quad \phi^*\cdot r\colon = \phi^*\lu s_R(r)\qquad 
\textrm{for}\ r\in R,\ \phi^*\in \asrd$$
and the right coaction
{
\be \tau_R:\asrd\to {\asrd_R}\ot {^R A}\qquad \phi^*\mapsto \sum_i 
\chi\inv(\pi_L \circ t_R \circ {}^*\kappa^j\circ S)\phi^* \ot k_j,\lb{taur}\ee}
where $\chi:\asrd\to \atld$ is the algebra anti-isomorphism (\ref{chi}). 
{
Note that the coactions (\ref{taul}) and (\ref{taur}) are independent of the
  choice of the dual bases.}

{
\begin{lem}
Let $\A=(\A_L,\A_R,S)$ be a Hopf algebroid and consider the algebra
anti-isomorphism $\chi$ in (\ref{chi}). For any $a,b\in A$ and $\psi_*\in
\atld$, 
\begin{eqnarray}
\lb{diff-key_2}
&&\chi\left(\chi^{-1}(\psi_*)\lu a\right)(b)=
\pi_L\circ t_R \circ \pi_R \left(s_L\circ \psi_*(a\di b)a\dii  \right)
\qquad \textrm{and}\\ 
\lb{diff-key}
&&\chi\left(\chi^{-1}(\psi_*)\lu S(a)\right)(b)=
\pi_L\left(a\ui t_L\circ \psi_*(S(a\uii)b)\right).
\end{eqnarray}
\end{lem}}

{
\pr
By the form of $\chi$ and its inverse, for $a,b\in A$ and $\psi_*\in \atld$, 
$$
\chi\left(\chi^{-1}(\psi_*)\lu a\right)(b)=
\pi_L\left(b\uii t_R \circ \pi_R (s_L\circ \psi_*(a\di {b\ui}\di)a\dii
   {b\ui}\dii) \right).
$$
For any $a\in A$ and $\psi_*\in \atld$, there is a well defined map
$A_L\ot {}_LA^R \ot{}^R A \to A$, $x\ot y\ot z \mapsto 
z t_R \circ \pi_R (s_L\circ \psi_*(a\di x) a\dii y)$.
Composing it with the equal maps $(\gamma_L \ot {}^R A)\circ \gamma_R = (A_L
\ot \gamma_R) \circ \gamma_L: A \to A_L\ot {}_LA^R \ot{}^R A$, we conclude
that 
$$
\chi\left(\chi^{-1}(\psi_*)\lu a\right)(b)=\pi_L\left({b\dii}\uii t_R \circ
\pi_R (s_L\circ \psi_*(a\di b\di)a\dii {b\dii}\ui) \right).
$$
Applying the right bialgebroid analogue of (\ref{2.6}) (in the first equality), 
counitality of $\gamma_R$ (in the second equality),
(\ref{2.5}), (\ref{hgdi}) and the left $R$-linearity of $\pi_R$ (in the third
equality), 
and (\ref{cros}) together with the counitality of $\gamma_L$ (in the last
equality), we conclude that 
\begin{eqnarray*}
\chi\left(\chi^{-1}(\psi_*)\lu a\right)(b)&=&
\pi_L\left({b\dii}\uii t_R \circ \pi_R (
t_R \circ \pi_R(s_L\circ \psi_*(a\di b\di) a\dii)
   {b\dii}\ui) \right)\\
&=&\pi_L\left(t_R \circ \pi_R(s_L\circ \psi_*(a\di
   b\di)a\dii) b\dii  \right)\\
&=&\pi_L\circ t_R \circ \pi_R \left(s_L\circ \psi_*(a\di
   b\di)a\dii s_L\circ \pi_L (b\dii)\right) \\
&=&\pi_L\circ t_R \circ \pi_R \left(s_L\circ \psi_*(a\di
   b)a\dii  \right). 
\end{eqnarray*}
Hence
$$
\chi\left(\chi^{-1}(\psi_*)\lu S(a)\right)(b)
=\pi_L\circ t_R \circ \pi_R \left(s_L\circ \psi_*(S(a)\di b)S(a)\dii \right)
=\pi_L\left (a\ui t_L\circ \psi_*(S(a\uii)b)\right),
$$
where the last equality follows by Proposition \ref{antip} (2).
\hfill \qed}

\bp \lb{lhm} 
Let $\A=(\A_L,\A_R,S)$ be a Hopf algebroid.

1)
Introduce the left $A$-module
$${_A \asrd}:\quad a\cdot \phi^*\colon =\phi^*\lu S(a)\qquad \textrm{for}\ a\in
A,\ \phi^*\in \asrd.$$
{
If the module $A_L$ is finitely generated and projective},
then $({_A\asrd},\tau_L)$ -- where $\tau_L$ is the map (\ref{taul}) -- is a
left-left Hopf module over $\A_L$.

2)
Introduce the right $A$-module
$${\asrd_A}:\quad \phi^*\cdot a\colon =\phi^*\lu a\qquad \textrm{for}\ a\in
A,\ \phi^*\in \asrd.$$
{
If the module ${}^R A$ is finitely generated and projective},
then $({\asrd_A},\tau_R)$ -- where $\tau_R$ is the map (\ref{taur}) -- is a
right-right Hopf module over $\A_R$.
\smallskip

The coinvariants of both Hopf modules $({_A\asrd},\tau_L)$ and
$(\asrd_A,\tau_R)$ are the elements
of $\il(\asrd)$.
\ep
\pr $\underline{1)}$: We have to show that $\tau_L$ is a left
$A$-module map. That is, for all $a\in A$ and $\phi^*\in \asrd$
\be \sum_i b_i\ot \chi\inv(\beta^i_*)\left(\phi^*\lu S(a)\right)=
\sum_i a\di b_i \ot \left(\chi\inv(\beta^i_*) \phi^*\right)\lu S(a\dii)
\lb{taumod} \ee
as elements of $A_L\ot {_L \asrd}$.
{
Note that for all $\phi^*,\psi^*\in \asrd$ and $a\in A$ 
\be (\phi^*\psi^*)\lu a=(\phi^*\lu a\uii) (\psi^*\lu a\ui). \lb{modalg}\ee 
Since $A_L$ is finitely generated and projective by assumption, using
(\ref{modalg}) and the dual basis property of $\{b_i\}$ and $\{\beta^i_*\}$,
one checks that (\ref{taumod}) is equivalent to the identity
$$
\left(\chi\inv(\psi_*)\lu s_L\ci \pi_L(a\di)\right)(\phi^*\lu S(a\dii))=
\sum_i \left(\chi\inv(\beta^i_*)\lu S\left[s_L\!\ci\!\psi_*(a\di
  b_i)\ a\dii\right]\right)(\phi^*\lu S(a_{(3)}))
$$
that is equivalent also to}
\be 
\chi\inv(\psi_*)\lu s_L\ci \pi_L(a)=
\sum_i \chi\inv(\beta^i_*)\lu S\left[s_L\!\ci\!\psi_*(a\di
  b_i)\ a\dii\right], \lb{diffi1}\ee
for all $a\in A$, $\psi_*\in \atld$.
{
Thus we have to prove (\ref{diffi1}). 
By (\ref{diff-key}), for any $x\in A$, 
\begin{eqnarray}
&&\sum_i\chi\left(
\chi^{-1}(\beta^i_*)\lu S(s_L\circ \psi_*(a\di b_i)a\dii)\right)(x)
\qquad\qquad\nn
&&\qquad\qquad\qquad\qquad
= \sum_i \pi_L\left(s_L\circ \psi_*(a\di b_i){a\dii}\ui t_L\circ \beta^i_* (
S({a\dii}\uii)x)\right) \nn
&&\qquad\qquad\qquad\qquad
= \sum_i \psi_*\left( t_L\circ \pi_L({a\dii}\ui t_L\circ \beta^i_* (
S({a\dii}\uii)x))a\di b_i\right).
\lb{L_3}
\end{eqnarray}
For any $x\in A$ and $\beta_* \in \atld$, there is a well defined map 
$A_L\ot {}_LA^R \ot {}^R A\to A$, $a\ot b \ot c\mapsto t_L\circ \pi_L (b
t_L \circ \beta_*(S(c)x))a$. Composing it with the equal maps 
$(A_L\ot \gamma_R)\circ \gamma_L= (\gamma_L \ot {}^R A)\circ \gamma_R:A \to
A_L\ot {}_LA^R \ot {}^R A$, we conclude that (\ref{L_3}) is equal to 
\begin{eqnarray*}
\sum_i \psi_*\left( t_L\circ \pi_L({a\ui}\dii t_L\circ \beta^i_* (
S(a\uii)x)){a\ui}\di b_i\right)
&=&\sum_i \psi_*\left(a\ui t_L\circ \beta^i_* (S(a\uii)x) b_i\right) \\
&=&\psi_*\left(a\ui S(a\uii)x \right)
=\psi_* \left(s_L\circ \pi_L(a) x\right).
\end{eqnarray*}
That is, $\sum_i \chi^{-1}(\beta^i_*)\lu S(s_L\circ \psi_*(a\di
b_i)a\dii) = \chi^{-1}(\psi_*(s_L\circ \pi_L(a)-))$.
Since for any $l\in L$, $\chi^{-1}\left(\psi_*(s_L(l)
-)\right)=\chi^{-1}(\psi_*)\lu s_L(l)$, this 
proves (\ref{diffi1}) hence claim 1).}

For the Hopf module $({_A \asrd},\tau_L)$, 
a projection onto the coinvariants is given by the
map (\ref{Em}), that takes the explicit form 
\be E_{\asrd}:\asrd\to \textrm{Coinv}(\asrd)\qquad \phi^*\mapsto
\sum_i\chi\inv(\beta^i_*)\phi^*\lu S^2(b_i). \lb{EA*}\ee
A left $s$-integral $\lambda^*$ on $\A$ is a coinvariant, since 
it is an invariant of the left regular $\asrd$-module and so for all
$a\in A$
$$ E_{\asrd}(\lambda^*)(a)=\sum_i
\chi\inv(\beta^i_*)(1_A)\lambda^*\left(S^2(b_i)a\right)=
\lambda^*\left[S^2\left(t_L\!\ci\!\beta^i_*(1_A)\ b_i\right) a
\right]= \lambda^*(a).$$
On the other hand, for all $a\in A$
\be  \sum_i S(b_i)(a\lu \beta^i_*)=S\left[t_L\!\ci\! \beta^i_*(a\di)\  
b_i\right]a\dii= s_R\ci\pi_R(a), \lb{auxi}\ee
hence for all $\phi^*\in \asrd$
\bea E_{\asrd}(\phi^*) \ru a&=&
\sum_i a\uii \ t_R\!\ci\! \pi_R\left\{\left[\phi^*\ru
S^2(b_i)a\ui\right]\lu\beta^i_*\right\}\nn 
&=& \sum_i t_R\!\ci\!\pi_R\!\ci\! S^2(b_i\uii)\ a\uii \ 
     t_R\!\ci\! \pi_R\left\{\left[\phi^*\ru S^2(b_i\ui)a\ui\right]\lu 
     \beta^i_*\right\}\nn
&=&\sum_i S(b_i\dii)\ \left(S^2(b_{i(1)})a\right)\uii \ 
    t_R\!\ci\! \pi_R\left\{\left[\phi^*\ru
    \left(S^2(b_{i(1)})a\right)\ui\right]\lu \beta^i_*\right\}\nn
&=&\sum_{i} S({b_i}\dii)\left\{\left[\phi^*\ru S^2({b_i}\di)a\right]\lu 
    \beta^i_*\right\}\nn
&=&\sum_{i,j} S(b_i)\left\{\left[\phi^*\ru S^2(b_j) a\right]\lu
  \beta^j_*\beta^i_*\right\}\nn 
&=&\sum_j s_R\ci\pi_R\left\{\left[\phi^*\ru S^2(b_j) a\right]\lu
  \beta^j_*\right\}
=s_R\ci E_{\asrd}(\phi^*)(a).
\nonumber\eea
That is, any coinvariant is an $s$-integral on $\A_R$. 
Here we used
(\ref{EA*}), the right analogue of (\ref{cros}), the identity $t_R\ci
\pi_R \ci S^2= S\ci s_R\ci \pi_R$, (\ref{hgdiv}), the right analogue
of (\ref{gmp}), the identity $\gamma_R\left[(\phi^*\ru a)\lu
\psi_*\right]=(\phi^*\ru a\ui)\lu \psi_*\ot a\uii$, holding true for
all $a\in A$, $\phi^*\in \asrd$ and $\psi_*\in \atld$, 
the dual basis property 
and (\ref{auxi}).

$\underline{2)}$:  
{
The proof is analogous to part 1), so we do not repeat the details.
We have to show that $\tau_R$ is a right $A$-module map. That is, for all
$a\in A$ and $\phi^*\in \asrd$ 
\be \sum_j \chi\inv(\pi_L \circ t_R \circ {}^*\kappa^j\circ S)\left(\phi^*\lu
a\right)\ot k_j= 
\sum_j \left(\chi\inv(\pi_L \circ t_R \circ {}^*\kappa^j\circ S)
\phi^*\right)\lu a\ui \ot k_j a\uii 
\lb{taurmod} \ee
as elements of ${\asrd_R} \ot {^R A}$. 
By similar steps used to show the equivalence of (\ref{taumod}) and
(\ref{diffi1}), the identity (\ref{taurmod}) is shown to be equivalent to 
\be
\chi\inv(\pi_L \circ t_R \circ {}^*\psi \circ S)  \lu t_R\circ \pi_R(a)=
\sum_j \chi\inv(\pi_L \circ t_R \circ {}^*\kappa^j\circ S)\lu a\ui s_R \circ
    {}^*\psi(k_j a\uii),
\lb{diffi2}
\ee
for all $a\in A$, ${}^*\psi\in \atrd$. Verification of (\ref{diffi2}) goes by
similar steps used to prove (\ref{diffi1}), making use of (\ref{diff-key_2}).
A projection onto the coinvariants is given by
$$
E:\phi^* \mapsto \sum_j \chi\inv(\pi_L \circ t_R \circ {}^*\kappa^j\circ
S)\phi^* \lu S(k_j). 
$$
By similar steps in part 1), one checks that $E(\lambda^*)=\lambda^*$, for
$\lambda ^*\in {\mathcal L}(\asrd)$, and $E(\phi^*)\in {\mathcal L}(\asrd)$,
for any $\phi^*\in \asrd$.}
\hfill\qed

Note that, if in a Hopf algebroid $\A=(\A_L,\A_R,S)$ both modules $A_L$ and
${}^R A$ are finitely generated and projective, then the dual bases 
$\{b_i\}\subset A$,
$\{\beta_*^i\}\subset \atld$ for $A_L$, and
$\{k_j\}\subset A$,
$\{{}^*\kappa^j\}\subset \atrd$ for ${}^R A$ are related via the identity
$$
\sum_i \beta_* ^i \ot S(b_i) = \sum_j \pi_L \circ t_R \circ {}^*\kappa_j
\circ S \ot k_j
$$
in $\atld^R \ot {}^R A$. In particular, in this case the projections $\asrd\to
{\mathcal L}(\asrd)$ in parts 1) and 2) of Proposition \ref{lhm} coincide and
(\ref{taur}) has the alternative form $\tau_R(\phi^*)=\sum_i
\chi^{-1}(\beta^i_*)\phi^* \ot S(b_i)$. 
\smallskip

Let us apply Theorem \ref{fundi} to the Hopf modules in Proposition \ref{lhm}.
If in a Hopf algebroid $\A=(\A_L,\A_R,S)$ the module $A_L$ is finitely
generated and projective, and the kernel of the map 
$$
{}^LA \ot \asrd^L\to {}^{\bf L} A \ot A_L\ot {}_L\asrd^{\bf L},\qquad 
a\ot \phi^* \mapsto a\ot \tau_L(\phi^*)- a\ot 1_A \ot \phi^*   
$$
is equal to ${}^L A \ot  {\mathcal L}(\asrd)^L$, then we conclude that 
\be\alpha_L:{^L A}\ot {\il(\asrd)^L}\to \asrd \qquad a\ot \lambda^*\mapsto
\lambda^*\lu S(a) \lb{alphal}\ee
is an isomorphism of left-left Hopf modules over $\A_L$. 
If the module ${}^RA$ is finitely generated and projective, and the kernel of
$$
{}^R \asrd \ot A_R \to {}^{\bf R}\asrd _R \ot {}^R A \ot A_{\bf R},\qquad 
\phi^* \ot a \mapsto \tau_R(\phi^*) \ot a - \phi^*\stac R 1_A \ot a
$$
is equal to ${}^R {\mathcal L}(\asrd) \ot A_R$, then we conclude that
\be \alpha_R:{^R\il(\asrd)}\ot A_R \to \asrd \qquad \lambda^*\ot a \mapsto
\lambda^*\lu a \lb{alphar}\ee
is an isomorphism of right-right Hopf modules over $\A_R$.
(The right $L$-module structure on $\il(\asrd)$ is given
by $\lambda^*\cdot l\colon =\lambda^*\lu s_L(l)$ and the left $R$-module
structure is given by $r\cdot \lambda^*\colon =\lambda^*\lu t_R(r)$ -- see the
explanation after (\ref{G}).)

\bc For a Hopf algebroid $\A=(\A_L,\A_R,S)$, such that all
of the modules $A^R$, ${^R A}$, ${_L A}$ and $A_L$ are 
finitely generated and projective, there exist non-zero elements in all of
$\il(\asrd)$, $\il(\atrd)$ $\ir(\asld)$ and $\ir(\atld)$.
\ec
\pr 
Since both modules $A^R$ and $A_L$ are finitely generated and projective
by assumption,
it follows from Proposition \ref{lhm} and Theorem
\ref{fundi} that the map (\ref{alphal}) is an isomorphism, hence there exist
non-zero elements in $\il(\asrd)$. 

For any element $\lambda^*$ of $\il(\asrd)$, $\lambda^*\ci S$ is a 
(possibly zero) element of $\il(\atrd)$ by Scholium \ref{dualsch}. Now
we claim that it is excluded by the bijectivity of the map
(\ref{alphal}) that $\lambda^*\ci S=0$ for all $\lambda^* \in
\il(\asrd)$. 
For if so, then by the surjectivity of the map
(\ref{alphal}) we have $\phi^*(1_A)=0$ for all $\phi^*\in \A^*$. But
this is impossible, since $\pi_R(1_A)=1_R$, by definition. 

It follows from Scholium \ref{dualsch}, 3.b) and 4.b)
that also $\ir(\asld)$ and $\ir(\atld)$ must contain non-zero
elements. 
\hfill\qed
\smallskip

Since none of the duals of a Hopf algebroid is known to be a Hopf
algebroid, it does not follow from Theorem \ref{fundi}, however,  that
for a Hopf algebroid, in which the total algebra is finitely generated and
projective as a module over the base algebra, also
$\il(A)$ and $\ir(A)$ contain non-zero elements. At the moment 
we do not know under what necessary conditions the existence of
non-zero integrals in a Hopf algebroid follows.

It is well known (\cite{Par}, Proposition 4) that the antipode of a finitely
generated and projective Hopf algebra over a commutative ring is bijective. We
do not know whether a result of the same strength holds true on Hopf
algebroids. Our present understanding of this question is formulated
in 
\bp \lb{bijantip}

For a Hopf algebroid $\A=(\A_L,\A_R,S)$, such that all of the modules $A^R$,
${^R A}$, ${_L A}$ and $A_L$ are finitely generated and projective, the
following statements are equivalent:

1) The antipode $S$ is bijective. 

2) There exists an invariant $\sum_k x_k\ot \lambda^*_k$ of the left $A$-module
   ${^R A}\ot {\il(\asrd)^R}$ -- defined via left multiplication in the first
   factor -- with respect to $\A_L$, satisfying in addition $\sum_k
   \lambda^*_k(x_k)=1_R$. (The right $R$-module 
   structure of $\il(\asrd)$ is defined by the restriction of the action
   on $(\asrd)^R$, i.e. as $\lambda^*\cdot r\colon =
   \lambda^*\left(\h t_R(r)\right)$.)
\ep
\pr
For any invariant $\sum_k x_k\ot \lambda^*_k$ of the left $A$-module
   ${^R A}\ot {\il(\asrd)^R}$  and any element $a\in A$, the identities
\bea \sum_k  S(a)x_k\ui \ot x_k\uii \ot \lambda^*_k&=&
\sum_k x_k\ui \ot a x_k\uii \ot \lambda^*_k\quad \textrm{and}\nn
\sum_k a x_k\ui \ot S(x_k\uii) \ot \lambda^*_k&=&
\sum_k x_k\ui \ot S(x_k\uii)a \ot \lambda^*_k
\nonumber\eea
hold true as identities in ${^R A^{\bf R}}\ot {^{\bf R} A}\ot {\il(\asrd)^R}$
and in ${^R A^{\bf R}}\ot {_{\bf R} A}\ot {\il(\asrd)^R}$, respectively.

$\underline{2)\Rightarrow 1)}$: In terms of the invariant $\sum_k x_k\ot
\lambda^*_k$, the inverse of the antipode is constructed explicitly as
$$ S\inv:A\to A\qquad a\mapsto \sum_k (\lambda^*_k\lu a)\ru x_k. $$

$\underline{1)\Rightarrow 2)}$ If $S$ is bijective then in the case of the
Hopf algebroid $\A_{cop}$ the isomorphism (\ref{alphal}) takes the form
$$\alpha_L^{cop}: A^L\ot {^L \il(\atrd)}\to \atrd\qquad a\ot {^*\lambda}\mapsto
 {^*\lambda}\ld S^{-1}(a),$$
where the left $L$-module structure on $\il(\atrd)$ is defined by $l\cdot
  {^*\lambda}\colon = {^*\lambda}\ld t_L(l)$.
 
In terms of $\sum_k x_k\ot {^*\lambda_k}\colon =(\alpha_L^{cop})\inv(\pi_R)$,
  the required invariant of ${^R A}\ot {\il(\asrd)^R}$ is given by $\sum_k
  x_k\ot  {^*\lambda_k}\ci S^{-1}$. 
\hfill\qed

\smallskip

In a Hopf algebroid $\A=(\A_L,\A_R,S)$, in which 
{
all of the modules $A_L$, ${}_LA$, $A^R$ and ${}^R A$ are}
finitely generated and projective, the  extensions $s_R:R\to A$ and
$t_L:L\op\to A$ satisfy the left depth two (or D2, for short) condition and the
extensions  $t_R:R\op\to A$ and $s_L:L\to A$ satisfy the right D2 condition of
\cite{KSz}. If furthermore  $S$ is bijective then all the four extensions
satisfy both the left and the right D2 conditions. This means
(\cite{KSz}, Lemma 3.7) in the case of $s_R:R\to A$, for example, the   
existence of finite sets (the so called D2 quasi-bases)
$\{d_k\}\subset A^R\ot {_R A}$, 
$\{\delta_k\}\subset {_R \End_R}({_R A^ R})$, 
$\{f_l\}\subset A^R\ot {_R A}$ and 
$\{\phi_l\}\subset  {_R\End_R}({_R A^ R})$ satisfying
\bea \sum_k d_k\cdot m_A\ci(\delta_k\ot {_R A})(u)&=&u\qquad \textrm{and}\nn
\sum_l m_A\ci(A^R\ot \phi_l)(u)\cdot f_l&=&u\nonumber\eea
for all elements $u$ in $A^R\ot {_R A}$, where the $A$-$A$ bimodule
structure on $A^R\ot {_R A}$ is defined by left multiplication in the
first factor and right multiplication in the second factor.

The D2 quasi-bases for the extension $s_R:R\to A$ can be constructed
in terms of the invariants $\sum_i x_i\ot \lambda^*_i\colon
=\alpha_L\inv(\pi_R)$ and $\sum_j x\pri_j\ot {^*\lambda\pri_j}\colon
=(\alpha_L^{cop})\inv (\pi_R)$ via the requirements that
\bea \sum_k d_k\ot \delta_k&=&\sum_i{x_{i(1)}}\ui\ot S({x_{i(1)}}\uii)\ot
[\lambda^*_i\lu S(x_{i(2)})]\ru \h \quad\qquad\textrm{and}\nn
\sum_l \phi_l\ot f_l&=&\sum_j \h\lu [x\pri_{j(1)}\ru \pi_L\ci s_R\ci
{^*\lambda\pri_j}\ci S\inv]\ot {x\pri_{j(2)}}^{(1)} \ot
S({x\pri_{j(2)}}^{(2)}) 
\nonumber\eea 
as elements of $A^R\ot {_R A^L}\ot {_L[{_R \End_R}({_R A^R})]}$ and of
$[{_R \End_R}({_R A^R})]_L\ot {_L A^R}\ot {_R A}$, respectively. (The
$L$-$L$ bimodule structure on ${_R \End_R}({_R A^R})$ is given by
$$ l_1\cdot \Psi \cdot l_2=s_L(l_1)\Psi(\h) s_L(l_2)\qquad
\textrm{for}\ l_1,l_2\in L,\ \Psi\in {_R \End_R}({_R A^R}).\qquad )$$
The D2 property of the extensions $t_R:R\op\to A$, $s_L:L\to A$ and
$t_L:L\op\to A$ follows by applying these formulae to the Hopf
algebroids $\A_{cop}$, $\A\op_{cop}$ and $\A\op$, respectively.

\smallskip
 
The following theorem and corollary, characterizing Frobenius Hopf algebroids 
$\A=(\A_L,\A_R,S)$ -- that is, Hopf algebroids  such that the
extensions, given by the source and target maps of the bialgebroids
$\A_L$ and $\A_R$, are Frobenius extensions --, are the main results of
this section. 

Recall that, for a homomorphism $s:R\to A$ of $k$-algebras, the canonical
$R$-$A$  bimodule ${_R A_A}$ is a 1-cell in the additive bicategory of
[$k$-algebras, bimodules, bimodule maps], possessing a right dual, the
bimodule ${_A A_R}$. If $A$ is 
finitely generated and projective as a left $R$-module, then 
${_R A_A}$ possesses also a left dual, the bimodule
${_A[{_R\Hom}({A},R)]_R}$ defined as
$$ a\cdot {\phi}\cdot r={\phi}(\h a)r\qquad \textrm{for}\ r\in R,\
a\in A,\ {\phi}\in {_R\Hom}(A,R). $$
A monomorphism of $k$-algebras $s:R\to A$ is called a {\em Frobenius
extension} if the module ${_R
A}$ is finitely generated and projective and the left and right duals
${_A A_R}$ and ${_A[{_R\Hom}(A,R)]_R}$ of the bimodule ${_R A_A}$ are
isomorphic. Equivalently, if $A_R$ is finitely generated and projective
and the left and right duals ${_R A_A}$ and ${_R [{\Hom_R}(A,R)]_A}$
of the bimodule ${_A A_R}$ are isomorphic.
This property holds if and only if there exists a {\em
Frobenius system} $(\psi,\sum_i u_i\ot v_i)$, where $\psi:A\to R$ is
an $R$-$R$ bimodule map and $\sum_i u_i\ot v_i$ is an element of
$A\stac{R} A$ such that
$$ \sum_i s\!\ci\!\psi(a u_i)\ v_i=a=\sum_i u_i\ s\!\ci\!\psi(v_i a)\qquad
\textrm{for\ all\ }a\in A.$$ 
\bt \lb{frobthm} 
If in a Hopf algebroid $\A=(\A_L,\A_R,S)$ all modules $A^R$, ${}^RA$, $A_L$ and
${}_L A$ are finitely generated and projective, then
the following statements are equivalent:

1.a) The map $s_R:R\to A$ is a Frobenius extension of $k$-algebras. 

1.b) The  map $t_R:R\op\to A$ is a Frobenius extension of $k$-algebras.
 
1.c) The map $s_L:L\to A$ is a Frobenius extension of $k$-algebras.

1.d) The  map $t_L:L\op\to A$ is a Frobenius extension of $k$-algebras.

2.a) The
     module ${\il(\asrd)}^L$, defined by $\lambda^*\cdot l\colon =
     \lambda^*\lu s_L(l)$, is free of rank 1.

2.b) $S$ is bijective and the 
     module ${^L \il(\atrd)}$, defined by $l\cdot
     {^*\lambda}\colon = {^*\lambda}\ld t_L(l)$, is free of rank 1.

2.c) The
     module ${_R \ir(\asld)}$, defined by $r\cdot {_*\rho}\colon =
     s_R(r) \rd {_*\rho}$, is free of rank 1.

2.d) $S$ is bijective and the 
     module ${\ir(\atld)_R}$, defined by 
     $\rho_*\cdot r\colon = t_R(r)\ru \rho_*$, is free of rank 1.

3.a) There exists an element $\lambda^*\in \il(\asrd)$ such that the map 
\be {\cal F}:A\to \asrd\qquad a\ \mapsto\  \lambda^*\lu a \lb{lambdal}\ee
     is bijective.

3.b) $S$ is bijective
     and there exists an element ${^* \lambda}\in \il(\atrd)$ such that the map 
$A\to \atrd,\quad a\ \mapsto\  {^* \lambda}\ld a $
     is bijective.

3.c) There exists  an element ${_*\rho}\in \ir(\asld)$ such that the map 
$A\to \asld,\quad a\ \mapsto\  a\rd {_*\rho} $
     is bijective.

3.d) $S$ is bijective
     and there exists an element ${\rho_*}\in \ir(\atld)$ such that the map 
$A\to \atld,\quad a\ \mapsto\  a\ru \rho_*  $
     is bijective.

4.a) There exists a left integral $\ell\in \il(\A)$ such that the map
\be {\cal F}^*:\asrd\to A\qquad \phi^*\ \mapsto\ \phi^*\ru \ell \lb{f*}\ee
is bijective.

4.b) $S$ is bijective and there exists a left integral $\ell\in
  \il(\A)$ such that the map 
\be {^*{\cal F}}: \atrd\to A\qquad {^*\phi}\ \mapsto\ {^*\phi}\rd \ell
  \lb{*F}\ee 
is bijective.

4.c) There exists a right integral $\err\in \ir(\A)$ such that the map
$ \asld\to A,\quad {_*\phi}\ \mapsto\ \err\ld {_*\phi}$
is bijective.

4.d) $S$ is bijective and there exists a right integral $\err\in
  \ir(\A)$ such that the map 
$ \atld\to A,\quad {\phi_*}\ \mapsto\ \err\lu {\phi_*}$
is bijective.

In particular, the integrals $\lambda^*$, ${^*\lambda}$, ${_*\rho}$ and
$\rho_*$ on $\A$ satisfying the
condition in 3.a), 3.b) 3.c) and 3.d), respectively, are Frobenius functionals
themselves for the extensions $s_R:R\to A$, $t_R:R\op\to A$, $s_L:L\to A$ and
$t_L:L\op\to A$, respectively.

What is more, under the equivalent conditions of the theorem the left
integrals $\ell\in \il(\A)$ satisfying the conditions 
in 4.a) and 4.b) can be chosen to be equal, that is, to be a
{\em non-degenerate} left integral in $\A$. Similarly, the right
integrals $\err\in \ir(\A)$ satisfying the conditions 
in 4.c) and 4.d) can be chosen to be equal, that is to be a
{\em non-degenerate} right integral in $\A$.
\et
\pr $\underline{4.a)\Rightarrow 1.a)}$: In terms of the left integral
$\ell$ in {\em 4.a)} define $\lambda^*\colon = {\cal F}^{*\ -1}(1_A)\in
\asrd$. 
The element $\ell\ot \lambda^*\in {^R\il(A)}\ot
{\il(\asrd)^R}$ is an invariant of the left $A$-module ${^R A}\ot
{\il(\asrd)^R}$, hence by Proposition \ref{bijantip} the antipode is
bijective. 
Since for all $\phi^*\in \asrd$
$$\phi^*\lambda^*={\cal F}^{*\ -1}(\phi^*\ru 1_A)=
{\cal F}^{*\ -1}(s^*\ci \phi^*(1_A)\ru 1_A)=
s^*\!\ci\! \pi^*(\phi^*)\ \lambda^*,$$
 $\lambda^*$ is
an $s$-integral on $\A_R$, so in particular an $R$-$R$ bimodule map
${_R A^R}\to R$. 

Since for all $a\in A$
$$ \ell\uii\ \  t_R\!\ci\! \lambda^*\left(S(a)\ell\ui\right)=a, $$
we have ${\cal F}^{*\ -1}(a)=\lambda^*\lu S(a)$.
{
Hence for all $\phi^*\in \asrd$ and $a\in A$,
$$\phi^*(a)=({\cal F}^{*\ -1}\circ {\cal F}^{*})(\phi^*)(a)=
\lambda^*\big(s_R\circ \phi^*(\ell\ui)S(\ell\uii)a\big)=
\phi^*\big(a\ell\ui s_R \circ \lambda^*\circ S(\ell\uii)\big).$$
Since $A^R$ is finitely generated and projective by assumption, this proves
that}
$\ell\ui\  s_R\!\ci\!\lambda^*\!\ci\! S(\ell\uii)=1_A$.
A Frobenius system for the extension $s_R:R\to A$ is
provided by $\left(\lambda^*,\ell\ui\ot S(\ell\uii)\right)$.

$\underline{1.a)\Rightarrow 2.a)}$: 
In terms of a Frobenius system
$\left(\psi, \sum_i u_i\ot v_i\right)$ for the extension $s_R:R\to A$, one
constructs an isomorphism of right $L$-modules as
\bea \kappa:& \il(\asrd)\to L\qquad \lambda^*&\mapsto\ \pi_L\left[\sum_i
s_R\!\ci\! \lambda^*(u_i)\  v_i\right] \lb{kappa}\\
\kappa^{-1}:&L\to \il(\asrd)\qquad l&\mapsto \ E_{\asrd}\left(\psi\lu
s_L(l)\right),
\eea

\vspace{-.73cm}
\noindent with inverse

\vspace{.2cm}
\noindent where $E_{\asrd}$ is the map (\ref{EA*}). 
The right $L$-linearity of $\kappa$ follows from the property of the
Frobenius system $(\psi, \sum_i u_i\ot v_i)$ that $\sum_i a u_i\ot
v_i=\sum_i u_i\ot v_i a$ for all $a\in A$, the
bialgebroid axiom (\ref{2.5}), and left $R$-linearity of the map
$\lambda^*:{_R A}\to R$ and the right $L$-linearity of $\pi_L:{_L A}\to L$.

The maps $\kappa$ and $\kappa\inv$ are mutual inverses as
\bea \kappa\inv \ci \kappa(\lambda^*)&=&
\sum_{i,j} [\chi\inv(\beta^j_*)\psi]\lu s_L\!\ci\!\pi_L
\left(s_R\!\ci\! \lambda^*(u_i)\ v_i\right)\ S^2(b_j)\nn
&=& \sum_{i,j} [\chi\inv(\beta^j_*)\psi]\lu S^2(b_j\uii)\
t_R\!\ci\!\pi_R \left[t_R\!\ci\! \pi_R\!\ci\!
S\left(s_R\!\ci\! \lambda^*(u_i)\ v_i\right)\ S^2(b_j\ui)\right]\nn
&=& \sum_{i,j} [\chi\inv(\beta^j_*)\psi]\lu S^2(b_j\uii)\
s_L\!\ci\pi_L\left[S(b_j\ui)\ s_R\!\ci\! \lambda^*(u_i)\ v_i\right]
=\lambda^*,\lb{kinvk}
\eea
where in the first step we used (\ref{modalg}), in the second step the
fact that by Proposition \ref{antip} we have $s_L\ci \pi_L= t_R\ci\pi_R\ci S$,
then the right analogue of (\ref{2.5}) and finally in the 
last step the identity in ${^R \il(\asrd)}\ot A_R$:
$$\sum_{i,j} [\chi\inv(\beta^j_*)\ \psi]\lu S^2(b_j\uii)\ot S(b_j\ui)\
s_R\!\ci\! \lambda^*(u_i)\ v_i =\alpha_R\inv\left(\sum_i \psi\lu s_R\!\ci\!
\lambda^*(u_i)\ v_i \right)=\lambda^*\ot 1_A,$$
which follows from the explicit form of the inverse of the map
(\ref{alphar}). In a similar way, also

\bea \kappa\ci\kappa\inv(l)&=&
\sum_{i,j} \pi_L\left[ s_R\!\ci\!\left(
\chi\inv(\beta^j_*)\ \psi\right)\left(s_L(l) S^2(b_j) u_i\right)\
v_i\right]\nn
&=& \sum_{i,j} \pi_L \left[ s_R\!\ci\!\left(
\chi\inv(\beta^j_*)\ \psi\right)\left(s_L(l) u_i\right)\
v_i\ S^2(b_j)\right]\nn
&=& \sum_{i,j} \pi_L \left[ s_R\!\ci\!\left(
\chi\inv(\beta^j_*)\ \psi\right)\left(s_L(l) u_i\right)\
v_i\ t_L\!\ci\!\pi_L \!\ci\!S^2(b_j)\right]\nn
&=& \sum_{i,j} \pi_L \left[ s_R\!\ci\!\left(
\chi\inv(\beta^j_*)\ \psi\right)\left(s_L(l)\ t_L\!\ci\!\pi_L
\!\ci\!S^2(b_j)\ u_i\right)\ v_i\right]\nn
&=& \sum_{i,j} \pi_L \left\{ s_R\!\ci\!\left[\left(
\chi\inv(\beta^j_*)\lu t_L\!\ci\!\pi_L
\!\ci\!S^2(b_j)\right)\psi\right](s_L(l) u_i)\
v_i\right\}= l,\nonumber\eea
where in the 
last step we used that $\sum_j  \chi\inv(\beta^j_*) \lu
t_L\!\ci\! \pi_L\!\ci\! S^2(b_j)=\chi\inv\left(\sum_j \beta^j_*\
t_*\!\ci\! \pi_L(b_j)\right)=\pi_R$. 

$\underline{2.a)\Rightarrow 3.a)}$: If $\kappa:\il(\asrd)^L\to L$ is an
isomorphism of $L$-modules then $\pi_R\ci s_L\ci \kappa: {^R \il(\asrd)}\to R$
is an isomorphism of $R$-modules.
Introduce the cyclic and separating generator
$\lambda^*\colon =\kappa\inv (1_L)$ for the module
$\il(\asrd)^L$. The map ${\cal F}$ in (\ref{lambdal}) is equal to
$\alpha_R\ci (\kappa\inv\!\ci \pi_L\!\ci\! t_R\  \ot\  
A_R)$ -- where $\alpha_R$ is the isomorphism (\ref{alphar}) -- hence it is
bijective. 

$\underline{3.a)\Rightarrow 4.a),4.b)}$: A Frobenius system for the extension
$s_R:R\to A$ is given in terms of the dual bases $\{b_i\}\subset A$ and
$\{\beta^*_i\}\subset \asrd$ for the module $A^R$ as $\left(\lambda^*,\sum_i
b_i\ot {\cal F}\inv(\beta^*_i)\right)$. 

The element $\ell\colon =\sum_i b_i\ \ t_L\!\ci\! \pi_L\!\ci\! {\cal
  F}\inv(\beta^*_i)$ is a left integral in $\A$. Using the identities
\bea \lambda^*\ru \ell&=& 
s_R\!\ci\! \lambda^*\left[\sum_i  b_i\ \  t_L\!\ci\! \pi_L\!\ci\! {\cal F}\inv
  (\beta^*_i)\right]=
t_L\!\ci\! \pi_L\left[\sum_i s_R\!\ci\! \lambda^*(b_i)\ \ 
{\cal F}\inv(\beta^*_i)\right]=1_A, 
\nn
\ell\ui\ot S(\ell\uii)&=& \sum_i b_i\ \ s_R\!\ci\! \lambda^*\left[{\cal
  F}\inv(\beta^*_i)\ell\ui\right]\ot S(\ell\uii)=
\sum_i b_i\ot S\left[\ell\uii \ \ t_R\!\ci\! \lambda^*(\ell\ui)\right]{\cal
  F}\inv(\beta^*_i)\nn
&=&\sum_i b_i\ot {\cal F}\inv(\beta^*_i),
\nonumber\eea
one checks that the inverse of the map ${\cal F}^*$ in (\ref{f*}) is given by
${\cal F}\ci S$.   
This implies, in particular, that $S$ is bijective.

The inverse of the map ${^*{\cal F}}$ in (\ref{*F}) -- defined in terms of the
same left integral $\ell$ -- is the map
$$A\to \atrd \qquad a\ \mapsto\ \lambda^*\ci S\ld S\inv(a).$$

$\underline{1.a)\Leftrightarrow 1.d)}$: The datum $\left(\psi,\sum_i u_i\ot
v_i\right)$ is a Frobenius system for the extension $s_R:R\to A$ if and only
if 
$\left(\pi_L\ci s_R\ci\psi, \sum_i u_i\ot v_i\right)$ is a Frobenius system 
for $t_L:L\op\to A$, where $\pi_L\ci s_R:R\to L\op$ was claimed to be an
isomorphism of $k$-algebras in part 1) of Proposition \ref{antip}. 

$\underline{1.a)\Rightarrow 1.c)}$: We have already seen that $1.a)\Rightarrow
3.a)\Rightarrow S$ is bijective. If the datum $\left(\psi,\sum_i u_i\ot
v_i\right)$ is a Frobenius system for the extension $s_R:R\to A$ then
$\left(\pi_L\ci s_R\ci \psi\ci S\inv,S(v_i)\ot S(u_i)\right)$ is a Frobenius
system for $s_L:L\to A$.

$\underline{4.c)\Rightarrow 1.c)\Rightarrow 2.c)\Rightarrow 3.c)\Rightarrow
  4.c)}$, $\underline{1.c)\Leftrightarrow 1.b)}$ and
$\underline{1.c)\Rightarrow 1.a)}$ follow by applying 
$4.a)\Rightarrow 1.a)\Rightarrow $ $2.a)\Rightarrow 3.a)\Rightarrow
  4.a)$, $1.a)\Leftrightarrow 1.d)$ and
${1.a)\Rightarrow 1.c)}$
to the Hopf algebroid $\A\op_{cop}$.

$\underline{1.b)\Rightarrow 2.b)\Rightarrow 3.b)\Rightarrow
  4.b)\Rightarrow 1.b)}$: We have seen that $1.b)\Leftrightarrow
1.c)\Rightarrow S$ is bijective. Hence we
can apply $1.a)\Rightarrow 2.a)\Rightarrow 3.a)\Rightarrow 
  4.a)\Rightarrow 1.a)$ to the Hopf algebroid $\A_{cop}$.

$\underline{1.d)\Rightarrow 2.d)\Rightarrow 3.d)\Rightarrow 
  4.d)\Rightarrow 1.d)}$ follows by applying ${1.b)\Rightarrow 2.b)\Rightarrow
    3.b)\Rightarrow 
  4.b)\Rightarrow 1.b)}$ to the Hopf algebroid $\A_{cop}\op$.
\hfill\qed
\smallskip

Based on Theorem \ref{frobthm}, we obtain the following generalization of
Theorem \ref{parthm}.

\bc\lb{cor:Frob}
For any Hopf algebroid $\A=(\A_L,\A_R,S)$, the following assertions are
equivalent. 

1.a) Both maps $s_R:R\to A$ and $t_R:R\op\to A$ are Frobenius extensions of
$k$-algebras. 

1.b) Both maps $s_L:L\to A$ and $t_L:L\op\to A$ are Frobenius extensions of
$k$-algebras. 

2.a) The module $A^R$ is finitely generated and projective and there exists an
element $\lambda^*\in \il(\asrd)$ such that the map  
$ {\cal F}:A\to \asrd,\quad a\ \mapsto\  \lambda^*\lu a $ is bijective.

2.b) $S$ is bijective, the module ${^R A}$ is finitely generated and projective 
     and there exists an element ${^* \lambda}\in \il(\atrd)$ such that the map 
$A\to \atrd,\quad a\ \mapsto\  {^* \lambda}\ld a $
     is bijective.

2.c) The module ${_L A}$ is finitely generated and projective and there exists
an element ${_*\rho}\in \ir(\asld)$ such that the map  
$A\to \asld,\quad a\ \mapsto\  a\rd {_*\rho} $
     is bijective.

2.d) $S$ is bijective, the module $A_L$ is finitely generated and projective 
     and there exists an element ${\rho_*}\in \ir(\atld)$ such that the map 
$A\to \atld,\quad a\ \mapsto\  a\ru \rho_*  $
     is bijective.

3.a) There exists a {\em non-degenerate} left integral, that is an element
$\ell\in \il(\A)$ such that both maps
${\cal F}^*:\asrd\to A,\quad \phi^*\ \mapsto\ \phi^*\ru \ell$ and 
${^*{\cal F}}: \atrd\to A,\quad {^*\phi}\ \mapsto\ {^*\phi}\rd \ell$
are bijective.

3.b) There exists a {\em non-degenerate right integral}, that is, an element
$\err\in \ir(\A)$ such that both maps 
$ \asld\to A,\quad {_*\phi}\ \mapsto\ \err\ld {_*\phi}$
and
$ \atld\to A,\quad {\phi_*}\ \mapsto\ \err\lu {\phi_*}$
are bijective.
\ec

\pr
$\underline{1.a)\Leftrightarrow 1.b)}$: This follows by the same reasoning
used to prove $1.a)\Leftrightarrow1.d)$ and $1.b)\Leftrightarrow1.c)$
in Theorem \ref{frobthm}. 

$\underline{1.a)\Rightarrow 2.a)}$: Since $s_R:R\to A$ is a Frobenius
extension by assumption, the modules $A^R$ and ${}_R A$ (hence also $A_L$) are
finitely generated and projective by definition. Similarly, since
$t_R:R^{op}\to A$ is a Frobenius extension, the modules ${}^R A$
and ${}_L A$ are finitely generated and projective. Thus this implication
follows by Theorem \ref{frobthm} $1.a)\Rightarrow3.a)$.

$\underline{2.a)\Rightarrow 3.a)\ and\ S\ is\ bijective}$: This is proven by
repeating the same steps used to prove $3.a)\Rightarrow4.a)$ and $4.b)$ in
Theorem \ref{frobthm}. 

$\underline{3.a)\Rightarrow 1.a)\ and\ S\ is\ bijective}$: Putting
$\lambda^*:= {\cal F}^{*-1}(1_A)$, the map $a\mapsto (\lambda^* \lu a)\ru
\ell$ is checked to be the inverse of $S$.

For any $r\in R$, ${\cal F}^*(r\lambda^*(-))=t_R(r)={\cal F}^*(\lambda^*\lu
s_R(r))$. So by the bijectivity of ${\cal F}^*$, we conclude that $\lambda^*$
is a left $R$-module map ${}_R A\to R$. Therefore the module ${}^R A$ is
finitely generated and projective with dual basis 
$ \lambda^*(S(-)\ell\ui) \ot \ell\uii \in \atrd_R\ot {}^RA$. 
The module $A_L$ is finitely generated and projective
by Lemma \ref{fgp}. Applying the same reasoning to the Hopf algebroid
$\A_{cop}$, we conclude that by the bijectivity of ${}^*{\cal F}$ also the
modules $A^R$ and ${}_L A$ are finitely generated and projective. Hence the
claim follows by Theorem \ref{frobthm}, $4.a)\Rightarrow1.a)$ and $1.b)$.

$\underline{1.b)\Leftrightarrow 2.c)\Leftrightarrow 3.b)}$: This follows by
applying $1.a)\Leftrightarrow 2.a)\Leftrightarrow 3.a)$ to the Hopf algebroid
$\A_{cop}^{op}$. 

$\underline{1.a)\Leftrightarrow 2.b)}$: Since we proved that $1.a)$ implies the 
bijectivity of the antipode, we can apply $1.a)\Leftrightarrow 2.a)$ to the
Hopf algebroid $\A_{cop}$. 

$\underline{1.b)\Leftrightarrow 2.d)}$: This follows by applying
$1.a)\Leftrightarrow 2.b)$ to the Hopf algebroid $\A_{cop}^{op}$.
\hfill\qed
\smallskip

It is proven in (\cite{HGD}, Theorem 5.17) that under the equivalent conditions
of Theorem \ref{frobthm} the duals, $\asrd$, $\atrd$, $\asld$ and $\atld$ of
the Hopf algebroid $\A$, possess (anti-) isomorphic Hopf algebroid structures.

The Hopf algebroids, satisfying the equivalent conditions of Theorem
\ref{frobthm}, provide examples of distributive Frobenius double algebras
\cite{Sz4}. (Notice that the integrals, which we call
non-degenerate, are called Frobenius integrals in \cite{Sz4}).  

Our result naturally raises the question, under what conditions on the base
algebra the equivalent conditions of Theorem \ref{frobthm} hold true. That
is, what is the generalization of Pareigis' condition -- the triviality of the
Picard group of the commutative base ring of a Hopf algebra -- to the
non-commutative base algebra of a Hopf algebroid. We are going to return to
this problem in a different publication.

\section{The Quasi-Frobenius property}
\lb{qf}
\setc{0}

It is known (\cite{Par}, Theorem added in proof), that any finitely
generated projective Hopf algebra over a commutative ring $k$ is 
(both a left and a right) quasi-Frobenius extension of $k$ in the sense of
\cite{Muller}.
In this section we examine in what Hopf algebroids is the total
algebra (a left or a right) quasi-Frobenius extension of the base
algebra. 

The quasi-Frobenius property of an extension $s:R\to A$ of $k$-algebras has
been 
introduced by M\"uller \cite{Muller} as a weakening of the Frobenius
property (see the paragraph preceeding Theorem \ref{frobthm}).
The extension $s:R\to A$ is {\em left quasi-Frobenius} (or left QF, for
short) if the module ${_R A}$ is finitely generated and projective 
(hence the bimodule ${_R A_A}$ possesses both a right dual ${_A A_R}$
and a left dual ${_A[{_R \Hom}(A,R)]_R}$ ) and the bimodule 
${_A A _R}$ is a direct summand in a finite direct sum of copies of
${_A[{_R \Hom}(A,R)]_R}$.

The extension $s:R\to A$ is {\em right QF} if $s$, considered as a map
$R\op\to A\op$, is a left QF extension. That is, if the module $A_R$ is
finitely generated and projective and the left dual bimodule ${_R
A_A}$ is  a direct summand in a finite direct sum of copies of the
right dual bimodule ${_R [{\Hom_R}(A,R)]_A}$. 

To our knowledge it is not known whether the notions of left and right
QF extensions are equivalent (except in particular cases, such as
central extensions, where the answer turns out to be affirmative
\cite{RosChas}; and Frobenius extensions, which are also both left and right
QF \cite{Muller}).

A powerful characterization of a Frobenius extension $s:R\to A$  is
the existence of a Frobenius system -- see the paragraph preceeding
Theorem \ref{frobthm}. 
In the following lemma a generalization to quasi-Frobenius extensions
is introduced:
\bl \lb{QFs} 
1) An algebra extension $s:R\to A$ is left QF if and only if the module ${_R A}$
is finitely generated and projective and there exist finite sets
$\{\psi_k\}$ $\subset {_R\Hom_R}(A,R)$ and $\{\sum_{i} u_i^k\ot v_i^k\}\subset
A\stac{R} A$ satisfying 
\bea &&\sum_{i,k}u_i^k \ s\!\ci\! \psi_k(v_i^k)=1_A\qquad\quad \textrm{and}\nn
&& \sum_{i} au_i^k\ot v_i^k=\sum_i u_i^k\ot v_i^k a\qquad
\textrm{for\ all\ values\ of\ } k \textrm{ and }a\in A. 
\nonumber\eea 
The datum $\{\psi_k,\sum_i u_i^k\ot v_i^k\}$ is called a {\em left
QF-system} for the extension $s:R\to A$.

2) An algebra extension $s:R\to A$ is right QF if and only if the module ${A_R}$
is finitely generated and projective and there exist finite sets
$\{\psi_k\}\subset {_R\Hom_R}(A,R)$ and $\{\sum_{i} u_i^k\ot v_i^k\}\subset
A\stac{R} A$ satisfying 
\bea &&\sum_{i,k}s\!\ci\! \psi_k(u_i^k)\ v_i^k=1_A\qquad\quad \textrm{and}\nn
&& \sum_{i} a u_i^k\ot v_i^k=\sum_i u_i^k\ot v_i^k a\qquad
\textrm{for\ all\ values\ of\ } k \textrm{ and }a\in A. 
\nonumber\eea 
The datum $\{\psi_k,\sum_i u_i^k\ot v_i^k\}$ is called a {\em right QF-system}
for the extension $s:R\to A$.
\el
\pr Let us spell out the proof in case $\underline{1)}$: Suppose that
there exists a left QF system 
$\{\psi_k,\sum_i u_i^k\ot v_i^k\}$ for the extension $s:R\to A$.
The bimodule ${_A A_R}$ is a direct summand in a finite direct sum of copies
of ${_A[{_R \Hom}(A,R)]_R}$ by the existence of $A$-$R$ bimodule maps
\bea \Phi_k:&{_R \Hom}(A,R)\to A\qquad&{\phi}\mapsto \sum_i
u_i^k\ s\!\ci\! {\phi}(v_i^k)\qquad \textrm{and}\nn
\Phi\pri_k:&A\to {_R \Hom}(A,R)\qquad &a\mapsto \psi_k(\h a)
\nonumber\eea
satisfying $\sum_k \Phi_k\ci \Phi_k\pri =A$.

Conversely, in terms of the $A$-$R$ bimodule maps $\{\Phi_k: 
{_R \Hom}(A,R)\to A\}$ and $\{\Phi\pri_k: A\to {_R \Hom}(A,R)\}$,
satisfying $\sum_k \Phi_k \ci \Phi_k\pri =A$,  
and the dual bases, $\{b_j\}\subset A$ and 
$\{{\beta_j}\}\subset {_R \Hom}(A,R)$ for the 
module ${_R A}$, a left QF system can be constructed as  
\bea &&\psi_k\colon = \Phi_k\pri(1_A)\in {_R \Hom_R}(A,R)\qquad
\qquad \textrm{and}\nn 
&& \sum_i u_i^k\ot v_i^k\colon = \sum_j \Phi_k({\beta_j})\ot b_j\in
A\stac{R} A.
\nonumber\eea

\vspace{-1cm}
\hfill\qed

\vspace{.2cm}

\smallskip
Lemma \ref{QFs} implies, in particular, that for a left/right QF extension $R\to
A$, $A$ is finitely generated and projective also as a right/left $R$-module.
\bt \lb{rightqf}
{
If in a Hopf algebroid $\A=(\A_L,\A_R,S)$ all modules $A^R$, ${}^R A$, $A_L$ and
${}_LA$ are finitely generated and projective, then the following}
are equivalent:

1.a) $s_R:R\to A$ is a left QF extension.

1.b) $t_L:L\op\to A$ is a left QF extension.

1.c) The module
  $\il(\asrd)^L$ -- defined by
  $\lambda^*\cdot l\colon = \lambda^*\lu s_L(l)$ -- is finitely
  generated and projective.

1.d) The module  $\il(\asrd)^L$ is flat.

1.e) The invariants of the left $A$-module ${^L A}\ot \il(\asrd)^L$ --
  defined via left multiplication in the first factor -- with respect
  to $\A_L$ are the elements of ${^L \il(\A)}\ot \il(\asrd)^L$.

1.f) There exist finite sets $\{ \ell_k\}\subset\il(\A)$ and
  $\{\lambda^*_k\}\subset \il(\asrd)$ satisfying $\sum_k
  \lambda^*_k\ci S(\ell_k)=1_R$.

1.g) The left $A$-module ${_A\asrd}$ -- defined by $a\cdot
  \phi^*\colon = \phi^*\lu S(a)$ -- is finitely generated and
  projective with generator set $\{\lambda^*_k\}\subset \il(\asrd)$. 

\noindent The following properties of $\A$ are also equivalent:

2.a) $s_L:L\to A$ is a right QF extension.

2.b) $t_R:R\op\to A$ is a right QF extension.

2.c) The module
  ${_R \ir(\asld)}$ -- defined by
  $r\cdot {_*\rho}\colon = s_R(r)\rd {_*\rho}$ -- is finitely
  generated and projective.

2.d) The module  ${_R \ir(\asld)}$ is flat.

2.e) The invariants of the right $A$-module ${_R \ir(\asld)}\ot A_R$ --
  defined via right multiplication in the second factor -- with respect
  to $\A_R$ are the elements of ${_R\ir(\asld)}\ot \ir(\A)_R$.

2.f) There exist finite sets $\{ \err_k\}\subset\ir(\A)$ and
  $\{{_*\rho _k}\}\subset \ir(\asld)$ satisfying $\sum_k
  {_*\rho _k}\ci S(\err_k)=1_L$.

2.g) The right $A$-module ${\asld_A}$ -- defined by ${_*\phi}\cdot
  a\colon = S(a)\rd {_* \phi}$ -- is finitely generated and
  projective with generator set $\{{_*\rho_k}\}\subset \ir(\asld)$. 

\noindent
If furthermore the antipode is bijective, then conditions
1.a)-1.g) and 2.a)-2.g) are equivalent to each other and also to 

1.h) The left $\atrd$-module on $A$ -- defined by ${^* \phi}\cdot a\colon
= {^* \phi}\rd a$ -- is finitely generated and projective with generator
set $\{\ell_k\}\in \il(\A)$.

2.h) The right $\atld$-module on $A$ -- defined by $a\cdot \phi_*\colon
= a\lu \phi_*$ -- is finitely generated and projective with generator
set $\{\err_k\}\in \ir(\A)$.
\et

\pr $\underline{1.a)\Leftrightarrow 1.b)}$: 
The datum $\{\psi_k,\sum_i u_i^k\ot v_i^k\}$ is a left QF system for the
extension $s_R:R\to A$ if and only if $\{\pi_L\ci s_R\ci \psi_k, \sum_i
u_i^k\ot v_i^k\}$ is a left QF system for $t_L:L\op\to A$.

$\underline{1.a)\Rightarrow 1.c)}$: 
In terms of the left QF system, $\{\psi_k,\sum_i u_i^k\ot v_i^k\}$ for 
the extension $s_R:R\to A$, the dual bases for the module ${
\il(\asrd)^L}$ are given with the help of the map (\ref{EA*}) as
$\{E_{\asrd}(\psi_k)\}\subset \il(\asrd)$ and 
$\{\kappa_k\colon = \pi_L\left[\sum_i s_R\!\ci\! \h (u_i^k)\
v_i^k\right]\}\subset {\Hom_L}({ \il(\asrd)^L},L)$. 

The right $L$-linearity of the maps $\kappa_k:\il(\asrd)\to L$ is checked
similarly to the right $L$-linearity of the map (\ref{kappa}). Notice that for
any $R$-$R$ bimodule map $\psi:{_R A ^R}\to R$ we have 
\bea 
E_{\asrd}(\psi\lu s_L(l))
&=&
\sum_j [\chi\inv(\beta^j_*)\psi]\lu s_L(l)S^2(b_j)\nn
&=&\sum_j [\chi\inv(t_*\!\ci\!\pi_L\!\ci\!t_R\!\ci\!\pi_R\!\ci\!t_L(l)
\ \beta^j_*)\psi]\lu S^2(b_j)\nn
&=&\sum_j [\chi\inv(\beta^j_*) \ t^*\!\ci\!\pi_R\!\ci\!t_L(l)\
\psi]\lu S^2(b_j)\nn
&=&\sum_j [\chi\inv(\beta^j_*) \ s^*\!\ci\!\pi_R\!\ci\!t_L(l)\
\psi]\lu S^2(b_j)\nn
&=&\sum_j [\chi\inv(s_*\!\ci\!\pi_L\!\ci\!t_R\!\ci\!\pi_R\!\ci\!t_L(l)
\ \beta^j_*)\ \psi]\lu S^2(b_j)\nn
&=& E_{\asrd}(\psi)\lu s_L(l),\nonumber
\eea
for all $l\in L$, where in the first step we used (\ref{EA*}) and
(\ref{modalg}), in the second step the property of the dual bases
$\{b_j\}\subset A$ and $\{\beta^j_*\}\subset \atld$ that $\sum_j
\beta^j_* \ot s_L(l) b_j = \sum_j t_*(l) \beta^j_*\ot b_j$ for all
$l\in L$ as elements of ${^L\atld}\ot A_L$, in the third step the
identity $\chi\inv \ci t_*=t^*\ci\pi_R\ci s_L$, in the fourth step the
fact that by the left $R$-linearity of $\psi$ we have $t^*(r)\psi=s^*(r)\psi$
for all $r\in R$, in the fifth step $\chi\inv \ci s_*=s^*\ci\pi_R\ci s_L$,
and finally $\sum_j \beta^j_* \ot b_j s_L(l) = \sum_j s_*(l)
\beta^j_*\ot b_j$, holding true for all 
$l\in L$ as an identity in ${^L\atld}\ot A_L$. 

The dual basis property of the sets $\{E_{\asrd}(\psi_k)\}$ and
$\{\kappa_k\}$ is verified by the property that
$\sum_{i,k} E_{\asrd}\left(\psi_k \lu s_L\ci
\kappa_k(\lambda^*)\right)=\lambda^*$ for all $\lambda^*\in
\il(\asrd)$, which is checked similarly to
(\ref{kinvk}). 

$\underline{1.c)\Rightarrow 1.d)}$ is a standard result.

$\underline{1.d)\Rightarrow 1.e)}$: 
{
Since the module $A_L$ is finitely generated and projective by assumption,}
the invariants of any left $A$-module
$M$ with  respect to $\A_L$ are the elements of the kernel of the map
$$ \zeta_M:M\to {^L\atld}\ot M_L\qquad m\mapsto \left(\sum_i
\beta^i_*\ot b_i\cdot m\right)-\pi_L\ot m,$$
where the right $L$
module $M_L$ is defined via $t_L$, and the sets $\{b_i\}\subset A$
and $\{\beta^i_*\}\subset \atld$ are dual bases for the module $A_L$.

The map $\zeta_A$, corresponding to the left regular $A$-module, is a
left $L$-module map ${^L A}\to {^{\bf L} \asrd}\ot {^L A_{\bf L}}$ and
$\zeta_{{^L A}\ot \il(\asrd)^L}=\zeta_A\ot \il(\asrd)^L$. Since
tensoring with $\il(\asrd)^L$ is an exact functor by assumption, it
preserves the kernels, that is the invariants in this case.

$\underline{1.e)\Rightarrow 1.f)}$: With the help of the map
(\ref{alphal}) introduce
$$\sum_k \ell_k\ot \lambda^*_k\colon = \alpha_L\inv(\pi_R)\in
\textrm{Inv}({^L A}\ot \il(\asrd)^L)\equiv {^L \il(\A)}\ot
\il(\asrd)^L.$$
It satisfies $\sum_k \lambda^*_k\ci S(\ell_k)=\alpha_L\ci
\alpha_L\inv(\pi_R)(1_A)=1_R$. 

$\underline{1.f)\Rightarrow 1.a)}$: In terms of the sets
$\{\ell_k\}\subset \il(\A)$ and $\{\lambda^*_k\}\subset \il(\asrd)$, a
left QF system for the extension $s_R:R\to A$ can be constructed as
$\{\lambda^*_k, \ell_k\ui\ot S(\ell_k\uii)\}$.

$\underline{1.f)\Rightarrow 1.g)}$: In terms of the sets
$\{\ell_k\}\subset \il(\A)$ and $\{\lambda^*_k\}\subset \il(\asrd)$, the
dual bases for the module ${_A \asrd}$ are given by
$\{\lambda^*_k\}\subset \il(\asrd)$ and $\{ \h \ru \ell_k\}\subset
{_A\Hom}({_A\asrd},A)$. 

$\underline{1.g)\Rightarrow 1.f)}$: In terms of the dual bases
$\{\lambda^*_k\}\subset \il(\asrd)$ and $\{\Xi_k\}\subset
{_A\Hom}({_A\asrd},A)$ one defines the required left integrals
$\ell_k\colon = \Xi_k(\pi_R)$ in $\A$.

The equivalence of the conditions $\underline{2.a)-2.g)}$ follows by
applying the above results to the Hopf algebroid $\A\op_{cop}$.

\noindent
Now assume that $S$ is bijective. Then

$\underline{1.f)\Leftrightarrow 2.f)}$ follows from Scholium
\ref{dualsch}.

$\underline{1.f)\Rightarrow 1.h)}$: Scholium \ref{intsch}, 1.b) and
Scholium \ref{dualsch}, 3.c)
can be used to show that in terms of the sets
$\{\ell_k\}\subset \il(\A)$ and $\{\lambda^*_k\}\subset \il(\asrd)$, the
dual bases for the left $\atrd$-module on $A$ are given by
$\{\ell_k\}\subset \il(\A)$ and $\{\lambda^*_k \ci S\ld S\inv(\h)
\}\subset {_{\atrd}\Hom}(A,\atrd)$. 

$\underline{1.h)\Rightarrow 1.f)}$: Let $\{\ell_k\} \subset \il(A)$
and $\{\chi_k\}\subset {_{\atrd}\Hom}(A,\atrd)$ be dual bases for the
left $\atrd$-module $A$. 
For any value of the index $k$, the element
$\chi_k(1_A)$ is an invariant of the left 
regular $\atrd$-module. Hence by the finitely generated projectivity of the
module ${}^R A$, it is
a $t$-integral on $\A_R$. By Scholium
\ref{dualsch} the elements $\lambda^*_k\colon =\chi_k(1_A)\ci S\inv$
are $s$-integrals on $\A_R$, satisfying
$$ \sum_k \lambda^*_k\ci S(\ell_k)=\pi_R[\sum_k \chi_k(1_A)\rd
\ell_k]=1_R.$$  

$\underline{2.f)\Leftrightarrow 2.h)}$ follows by applying
$1.f)\Leftrightarrow 1.h)$ to the Hopf algebroid $\A\op_{cop}$.
\hfill\qed

\smallskip

If the antipode of a Hopf algebroid $\A=(\A_L,\A_R,S)$ is bijective
then the application of Theorem \ref{rightqf} to the Hopf algebroid
$\A\op$ results equivalent conditions under which the extensions
$s_R:R\to A$ and $t_L:L\op\to A$ are right QF, and $s_L:L\to A$ and
$t_R:R\op\to A$ are left QF.

\smallskip

In order to show that -- in contrast to Hopf algebras over commutative
rings -- not any finitely generated projective Hopf algebroid is
quasi-Frobenius, let stand here an example (with bijective antipode)
such that the total algebra is finitely generated and projective as a
module over the base algebra (in all the four senses listed in
(\ref{amod}) ) and the total algebra is neither a left nor a right QF
extension of the base algebra.

The example is taken from (\cite{Lu}, Example 3.1) where it is shown
that for any algebra $B$ over a commutative ring $k$ the $k$-algebra
$A\colon =B\stac{k} B\op$ has a left bialgebroid structure denoted by $\A_L$,
over the base $B$, with structural maps

\begin{center}
\begin{tabular}{lll}
$s_L:$&$B\to A\qquad $&$b\mapsto b\ot 1_B$\\
$t_L:$&$B\op \to A\qquad $&$b\mapsto 1_B\ot b$\\
$\gamma_L:$&$ A\to A_B \ot {_B A} \qquad$
&$b_1\ot b_2\mapsto (b_1\ot 1_B)\ot (1_B\ot b_2)$\\
$\pi_L:$&$A\to B \qquad $&$b_1\ot b_2\mapsto b_1b_2.$
\end{tabular}
\end{center}
\vspace{-.7cm}
\be\lb{Lul}\ee
The bialgebroid $\A_L$ satisfies the
Hopf algebroid axioms of \cite{Lu} with the involutive antipode $S$,
equal to the flip map 
\be S:B\stac{k} B\op\to B\op\stac{k} B\qquad b_1\ot b_2\mapsto b_2\ot
b_1. \lb{LuS}\ee
The reader may check that $A$ has a Hopf algebroid structure also in
the sense of this paper with left bialgebroid structure
(\ref{Lul}),  antipode (\ref{LuS}) and right bialgebroid structure
$\A_R=(A,B\op, S\ci s_L,S\ci t_L, (S\ot S)\ci \gamma_L\op\ci S, \pi_L\ci
S)$.

If $B$ is finitely generated and projective as a $k$-module, then all
modules $A^{B\op}$, ${^{B\op} A}$, $A_B$ and ${_B A}$ are finitely
generated and projective. What is more, we have
\bl \lb{ex}
Let $B$ be an algebra over a commutative ring $k$ with trivial
center. The following statements are equivalent:

1) The extension $k\to B$ is left QF.

2) The extension $k\to B$ is right QF.

3) The extension $B\to B\stac{k} B\op,\quad b\mapsto b\ot 1_B$  is
   left QF.

4) The extension $B\to B\stac{k} B\op,\quad b\mapsto b\ot 1_B$ is
   right QF. 
\el
The equivalence $\underline{1)\Leftrightarrow 2)}$ is proven in
\cite{RosChas} and  the rest can be proven using the technics
of quasi-Frobenius systems.

In view of Lemma \ref{ex}, it is easy to construct a finitely
generated projective Hopf algebroid which is not QF. Let us choose, for
example, $B$ to be the  algebra of $n\times n$ upper triangle
matrices with entries in a commutative ring $k$. Then $B$ has a
trivial center and it is neither a left nor a right QF extension of $k$.
Hence $A=B\stac{k} B\op$ is neither a left nor a right QF extension of $B$.


\end{document}